\documentclass[english,10pt,a4paper]{article}
\usepackage{amsmath,amssymb}
\usepackage{amssymb,cite}
\usepackage{graphicx}
\usepackage[hidelinks]{hyperref}
\usepackage{tikz}
\usepackage{epstopdf}
\usepackage[total={6.5in, 10in}]{geometry}

\definecolor{myblue}{RGB}{55,126,184}
\definecolor{myred}{RGB}{228,26,28}
\definecolor{mygreen}{RGB}{77,175,74} 

\definecolor{col1}{RGB}{127,205,187}
\definecolor{col2}{RGB}{65,182,196}
\definecolor{col3}{RGB}{29,145,192}
\definecolor{col4}{RGB}{34,94,168}
\definecolor{col5}{RGB}{37,52,148}
\definecolor{col6}{RGB}{8,29,88}

\title{Isolated resonances and nonlinear damping}
\author{Giuseppe Habib$^{1}$, Giuseppe I. Cirillo$^{2}$, Gaetan Kerschen$^{3}$}
\begin{document}
\date{}
\maketitle
{\em $^{1}$Department of Applied Mechanics, Budapest University of Technology and Economics, Hungary.}

{\em $^{2}$Department of Engineering, University of Cambridge, UK.}

{\em $^{3}$Department of Aerospace and Mechanical Engineering, University of Liege, Belgium.}

\bigskip Author preprint version -- manuscript under review"


\abstract{We analyze isolated resonance curves (IRCs) in a single-degree-of-freedom system with nonlinear damping. The adopted procedure exploits singularity theory in conjunction with the harmonic balance method. The analysis unveils a geometrical connection between the topology of the damping force and IRCs. Specifically, we demonstrate that extremas and zeros of the damping force correspond to the appearance and merging of IRCs.}

\section{Introduction}

Isolated resonance curves (IRCs) are branches of periodic solutions that are disconnected from the main resonance branch in the frequency response of nonlinear oscillating systems.
One of the first studies addressing their existence in engineering systems appeared already in 1955 \cite{abramson1955}, when Abramson illustrated the existence of a detached branch in the frequency response of a softening Duffing oscillator.
In 1966, Bouc \cite{bouc1964influence} showed that the resonance curve of an electrical resonator can present isolated portions if the iron core has hysteresis due to saturation.
The same system was studied two years later by Hayashi \cite{hayashi1966influence} adopting a harmonic balance approach limited to a single frequency. His analytical investigation allowed for a physical interpretation of the phenomenon.
In 1970, Hagedorn \cite{hagedorn1970parametric} identified detached resonant curves in a parametrically excited oscillator encompassing nonlinear damping and elastic force.
From 1968 to 1973, Iwan and Furuike \cite{iwan1968steady,furuike1971dynamic, iwan1973transient} studied the effect of hysteresis in mechanical systems in relation to IRCs. In particular, they analyzed the dynamics of a harmonically excited single-degree-of-freedom (DoF) system with a limited slip joint. The appearance of disconnected response curves was predicted analytically adopting a method based on the dissipation function. The results demonstrated that, for sufficiently small linear damping, IRCs are present even for very small slip. Furthermore, an estimation of the basins of attraction of the IRCs was provided, illustrating their relatively high robustness and, therefore, practical relevance.
In 1975, Koenigsberg and Dunn \cite{koenigsberg1975jump} rediscovered IRCs in electrical systems, while studying the jump phenomenon in a single-DoF inertial gyro employing ternary rebalance logic, where a nonlinearity was given by a relay with a deadband. ``Jump resonant frequency islands" were analytically and experimentally (by analog computer) identified. In light of the existence of IRCs, described in \cite{koenigsberg1975jump} as a new phenomenon, and in order to gain further insight into the dynamics, the employment of amplitude-frequency-excitation surfaces, instead of the classical amplitude-frequency curves, was suggested.
In 1978, Hirai and Sawai \cite{hrai1977jump,hirai1978general} established the first classification of IRCs. Adopting a graphical approach, a methodology to predict the existence of IRCs in nonlinear systems was developed. The procedure allows one to distinguish between inner and outer IRCs (laying below or above the main resonance curve branch), which were called lakes and islands, respectively.
This classification categorizes IRCs as simple islands, island chains, over-and-under-islands, boomerang islands, double island and islands with a lake.
Analytical results were confirmed by experiments on an analog computer for a system with deadband.
In the same year, Fukuma and Matsubara \cite{fukuma1978jump} developed a procedure to predict different kinds of IRCs based on the concept of resonance response surface. Very different types of nonlinearities were considered and involved scenarios, including IRCs, were illustrated for models of a steam-heat exchanger and of a pipeline system with recycle.

The aforementioned papers provide a solid basis for the analysis of IRCs in nonlinear systems.
Despite their fundamental contributions, these papers have received very little attention in the technical literature, and one objective of the present paper is to bring them back to light. Since these seminal studies, very few papers were published on the topic until the beginning of the 21th century, which has witnessed a resurgence of interest in IRCs.
Capecchi and Vestroni \cite{capecchi1990periodic} encountered an IRC investigating the dynamics of a single-DoF system with hysteretic elastic force.
IRCs were disclosed for a bilinear model of a suspension bridge \cite{doole1996piece}.
A family of sub-harmonic IRCs was also found in a piecewise linear single-DoF system \cite{duan2007sub}, where discontinuity seemed to be a critical parameter for their generation \cite{duan2008isolated}.
Elmeg{\aa}rd et al. \cite{elmegaard2014bifurcation} and Bureau et al. \cite{bureau2014} numerically and experimentally illustrated the generation of IRCs caused by impact and subharmonic resonances.
Atomic force microscopies in tapping-mode operation can also undergo IRCs related to a discontinuity \cite{lee2003nonlinear, misra2010degenerate}.
Nayfeh and Mook \cite{nayfeh1995nonlinear} illustrated the presence of IRCs related to the interaction of subharmonic and superharmonic resonances in a two-DoF system.
Perret-Liaudet and Rigaud \cite{perret2007superharmonic} showed the existence of IRCs in a softening single-DoF system in correspondence of a superharmonic resonance.
Rega \cite{rega2004nonlinear} denoted an elongated IRC in the vicinity of the 1/3 subharmonic region of a suspended cable.
Lenci and Ruzziconi \cite{lenci2009nonlinear} showed the appearance of an IRC far from any resonance, while studying a single-DoF model of a suspended bridge, which includes quadratic and cubic hardening nonlinearity.
DiBerardino and Dankowicz \cite{diberardino2014accounting} identified IRCs related to symmetry breaking in a two-DoF system; adopting isola singularity identification, in combination with a multiple scale approach, they managed to predict the outbreak of an IRC.
In \cite{arroyo2016duffing}, an IRC related to an internal resonance was encountered, while in \cite{mangussi2016internal} the phenomenon was explained through a frequency gap due to phase locking.

Isolated branches of periodic solutions can also be related to self-excited oscillations  \cite{takacs2008isolated, luongo2011parametric, zulli2012bifurcation, dimitriadis2017introduction}. Apart from mechanical systems, this is a typical and well-studied phenomenon in chemical reactors \cite{van1953autothermic, hlavavcek1970modeling, 
uppal1976classification, razon1987multiplicities}, but also in biological models \cite{doedel1984computer, pavlou1992microbial}.
However, in these cases they are not necessarily related to a resonance between external excitation and system response, therefore the phenomenon is qualitatively different from the one considered here.

In the last 15 years numerous studies about nonlinear vibration absorbers appeared.
Although very different types of vibration absorbers exist, most of them exploit internal resonances, making them prone to generation of IRCs.
The nonlinear energy sink (NES), consisting of a purely nonlinear resonator, if attached to single- or multi-DoF primary systems, presents IRCs \cite{starosvetsky2008response, starosvetsky2008dynamics}, whose existence was verified also experimentally \cite{gourc2014experimental}. An attempt to eliminate this undesired phenomenon demonstrated that IRCs can be avoided if the absorber has a properly-tuned piecewise-quadratic damping characteristic \cite{starosvetsky2009vibration}.
Similarly, the nonlinear tuned vibration absorber (NLTVA), possessing both a linear and a nonlinear elastic force characteristic, can present an IRC that limits its range of operation \cite{habib2015nonlinear, habib2016principle}.
A numerical procedure exploiting bifurcation tracking allowed to define regions of appearance of the IRC \cite{detroux2015performance}.
Alexander and Schilder \cite{alexander2009exploring} demonstrated that the existence of an IRC completely jeopardizes the efficiency of the NLTVA if applied to a linear primary system.
Cirillo et al. \cite{cirillo2016analysis}, analyzing the system's singularities, managed to eliminate the detrimental IRC using an NLTVA that includes higher order nonlinearities.

Partially motivated by the existence of IRCs in vibration absorbers, several studies dedicated to their identification and prediction were published in the last decade.
Gatti et al. \cite{gatti2010interaction, gatti2010response,  gatti2011effects} thoroughly studied the generation of IRCs in a two-DoF system, consisting of a linear primary system with a nonlinear attachment (resembling the NES). Due to the very small mass ratio, the system is subject to a sort of ideal base excitation, which explains the generation of IRCs when it undergoes a 1:1 resonance.
IRCs were predicted studying the parameter values for which the system response is multivalued.
In \cite{gatti2016uncovering} the analysis was extended to the case of large mass ratio, while in \cite{gatti2017inner} experimental proof of the existence of inner IRCs was given.
IRCs related to internal resonances were predicted through bifurcation tracking in \cite{detroux2015harmonic}. The result of the continuation of the fold bifurcations in the amplitude-frequency-excitation space resembles the amplitude-frequency-excitation surface proposed 40 years earlier by Koenigsberg and Dunn \cite{koenigsberg1975jump}.
Kuether et al. \cite{kuether2015nonlinear} developed a numerical technique based on nonlinear normal modes (NNMs) and energy balance to predict the occurrence of IRCs related to internal resonances; their procedure was tested on a cantilever beam.
Hill et al. \cite{hill2016analytical} further developed this method, allowing for a fully analytical approach.
It was then implemented on a two-DoF nonlinear system \cite{hill2016analytical},  on a cantilever beam with a nonlinear spring at its free end \cite{shaw2016periodic} and on a model of a pinned-pinned beam \cite{hill2017identifying}.
In \cite{noel2015isolated, detroux2016IMAC}, IRCs related to a 3:1 resonance of a two-DoF system were studied numerically and experimentally.


Most of the cited papers agree that the identification of IRCs is particularly troublesome, since standard continuation techniques tend to overlook them. Similarly, experimental or numerical investigations performed through frequency sweep are also unable to catch them if not carried out thoroughly.
An illustrative example is given in \cite{gourdon2007nonlinear}, where it is stated that the NES achieves better performance than a linear vibration absorber because of an undetected IRC \cite{alexander2009exploring}.

An alternative strategy to reveal IRCs is provided by singularity theory, which offers a very well-established mathematical framework to investigate the appearance and disappearance of IRCs.
In fact, the onset of an IRC corresponds to an isola singularity in the frequency response function, while its merging with another branch corresponds to a simple bifurcation (see for example \cite{spence1984numerical, golubitsky1985, troger1991nonlinear, drazin1992nonlinear, janovsky1992computer, seydel2009practical}).
In spite of this, only few researchers \cite{diberardino2014accounting, cirillo2016analysis} implemented singularity analysis for IRC identification.

In this study we attempt to fully exploit the potential of singularity theory for the prediction and identification of IRCs.
In particular, we analyze their relation with the topology of the nonlinear damping force, which suggests the presence of general rules for their appearance.
The system considered has a single DoF with a damping force given by a smooth symmetric function.
The simplicity of the system allows us to isolate the effect of damping, eliminating any contribution due to modal interaction, discontinuity, friction or symmetry breaking, which are known to generate IRCs.

\section{Numerical evidence of an isolated resonance}

We consider a harmonically-excited single-DoF oscillator, possessing a nonlinear damping force. Its dynamics is governed by the differential equation
\begin{equation}\label{eq:original_damp}
  m y''+k y+\tilde c_1 y'+\tilde c_3 y^{\prime3}+\tilde c_5 y^{\prime5}=2\tilde f \cos\left(\tilde\omega \tilde t\right),
\end{equation}
which, applying a standard nondimensionalization procedure, can be reduced to 
\begin{equation}\label{NL_DF}
\ddot x+x+c_1\dot x+c_3\dot x^3+\dot x^5=2 f\cos\left(\omega t\right),
\end{equation}
where  $\omega=\tilde\omega\sqrt{m/k}$, $t=\tilde t\sqrt{k/m}$, $c_1=\tilde c_1/\sqrt{km}$, $c_3=\tilde c_3k^{-1/4}m^{-1/4}\tilde c_5^{-1/2}$, $x=\tilde c_5^{1/4}k^{3/8}m^{-5/8}y$ and $f=\tilde c_5^{1/4}m^{-5/8}k^{-5/8}\tilde f$.
Linear and quintic terms are assumed strictly positive ($\tilde c_1>0$, $k>0$ and $\tilde c_5>0$), such that the trivial positions of the system remain locally stable and the trajectories remain bounded. Furthermore $m>0$.
Dots indicate derivation with respect to the dimensionless time $t$.

\begin{figure}
\begin{centering}
\setlength{\unitlength}{\textwidth}
\begin{picture}(1,0.35)
\put(0.01,-0.01){\includegraphics[trim = 12mm 10mm 15mm 10mm,clip,width=0.32\textwidth]{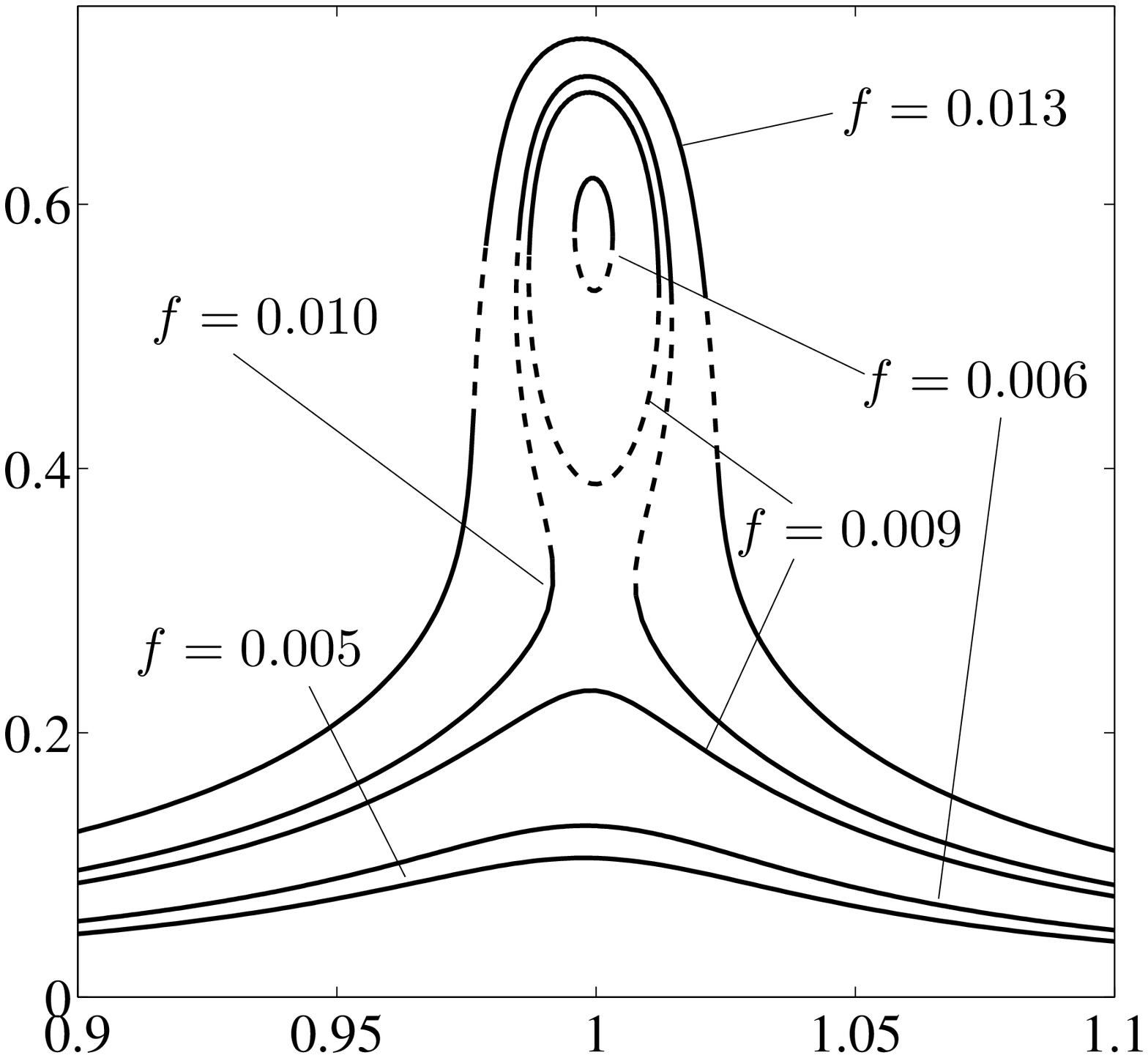}}
\put(0.35,-0.01){\includegraphics[trim = 09mm 10mm 18mm 10mm,clip,width=0.32\textwidth]{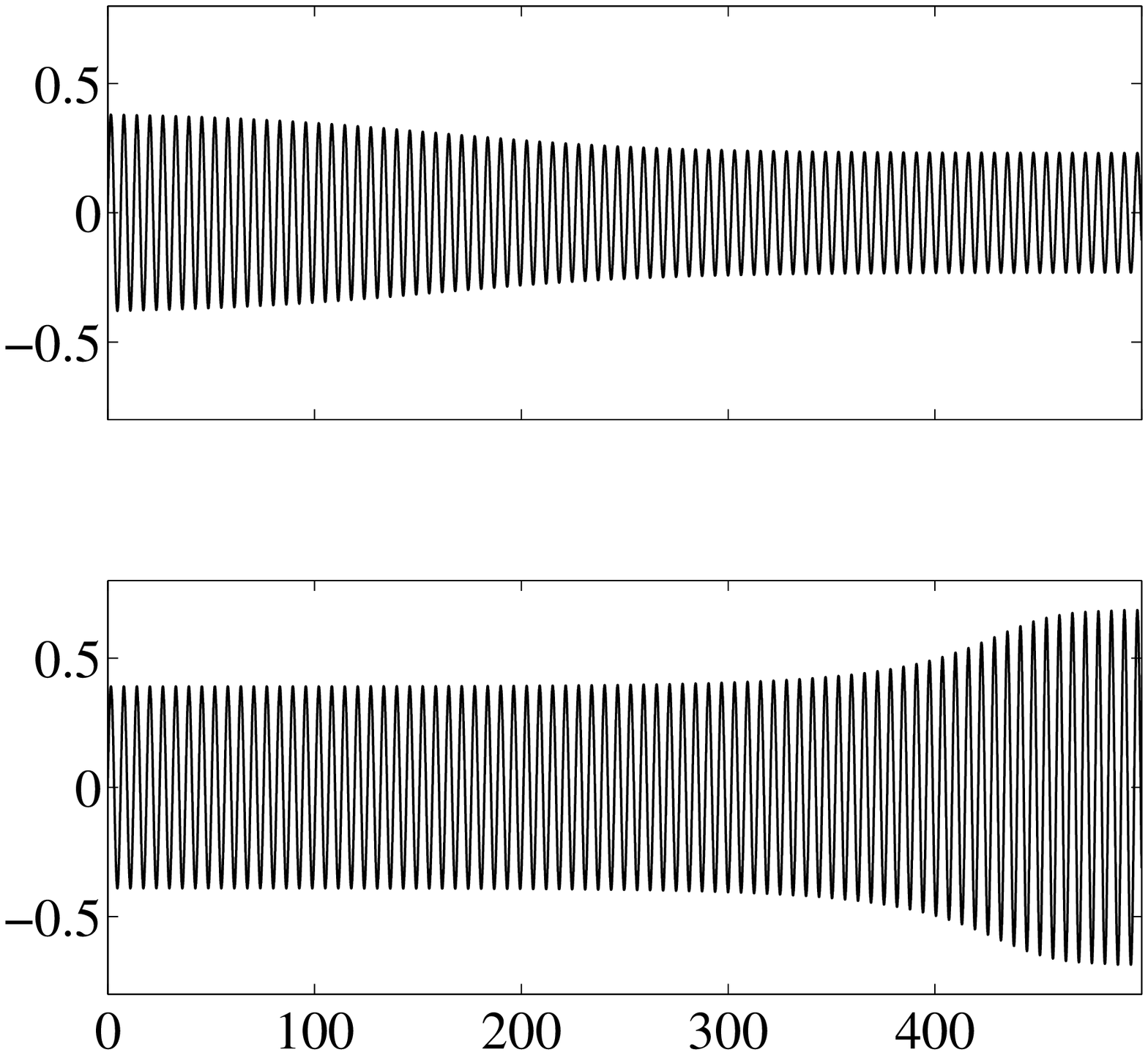}}
\put(0.69,-0.01){\includegraphics[trim = 12mm 10mm 15mm 10mm,clip,width=0.32\textwidth]{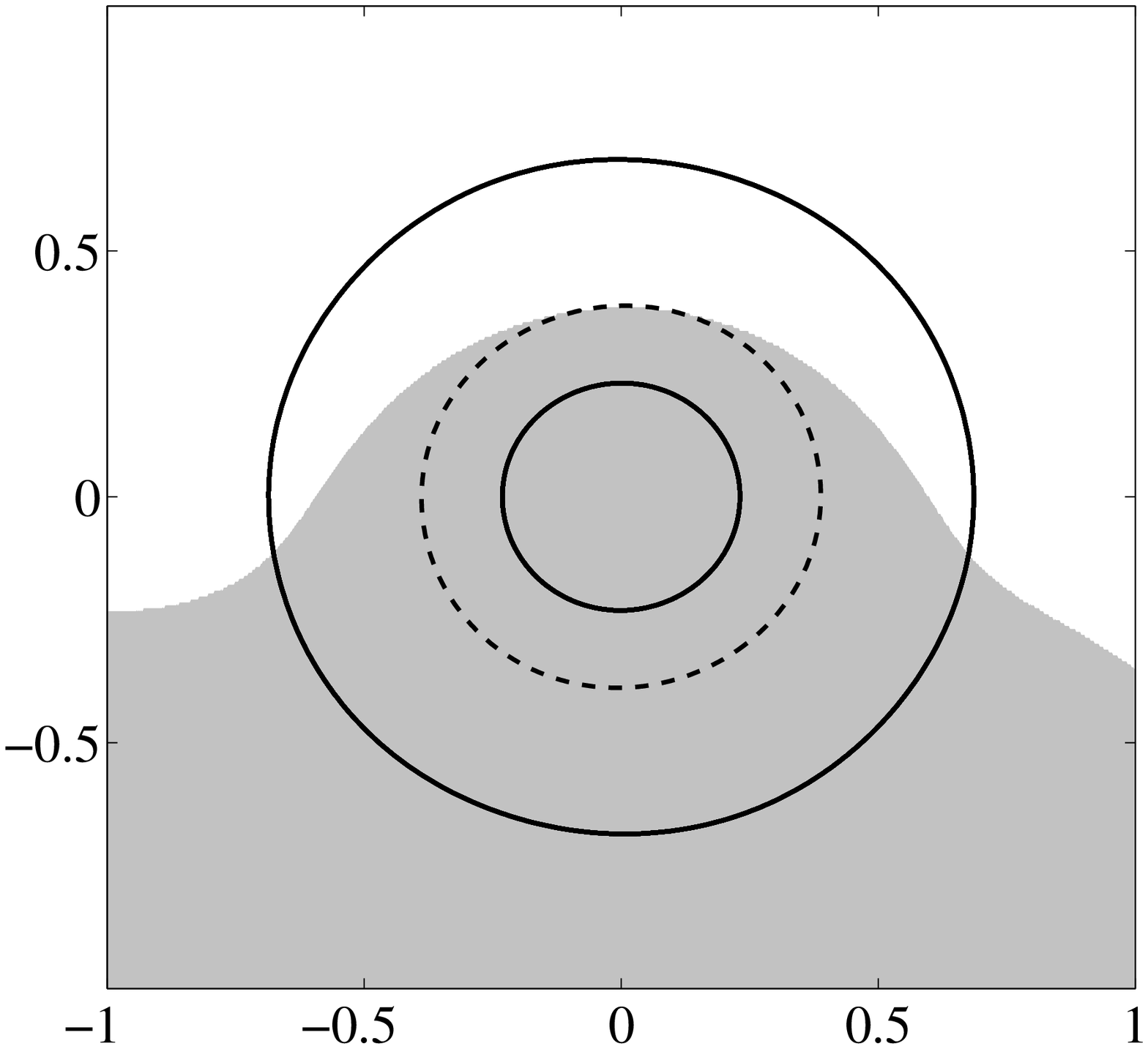}}
\put(0.04,0.26){{\textbf (a)}}
\put(0.39,0.26){{\textbf (b)}}
\put(0.39,0.103){{\textbf (c)}}
\put(0.72,0.26){{\textbf (d)}}
\put(0.175,-0.03){$\omega$}
\put(0.51,-0.03){$t$}
\put(0.855,-0.03){$x$}
\put(-0.01,0.14){\rotatebox{90}{$x$}}
\put(0.34,0.223){\rotatebox{90}{$x$}}
\put(0.34,0.061){\rotatebox{90}{$x$}}
\put(0.68,0.14){\rotatebox{90}{$\dot x$}}
\end{picture}
\end{centering}
\caption{(a) Frequency response for different values of the forcing amplitude as indicated in the figure; (b,c) time series for initial conditions $x(0)=0$, $\dot x(0)=0.38$ (b) and $\dot x(0)=0.39$ (c); (d) existing periodic solutions (solid lines: stable, dashed lines unstable) and basin of attraction of the two stable motions (shaded area corresponds to the smaller attractor). For the numerical computations $c_1=0.1$, $c_3=-0.6$, $f=0.009$ and $\omega=1$.}\label{numerical_results}
\end{figure}

Adopting a combination of direct numerical simulations and continuation techniques at different forcing amplitudes, it can be verified that, for certain parameter values, an IRC exists.
A typical scenario that can be obtained is as follows.
For low forcing amplitude the system behaves similarly to a linear one and the resonance curve presents the classical bell shape (Fig.~\ref{numerical_results}a for $f=0.005$).
Increasing the forcing amplitude the IRC appears, in the considered case this occurs for $f$ between 0.005 and 0.006.
A further increase in the forcing amplitude causes an enlargement of the IRC, which gets closer and closer to the main branch of the resonance curve, until they merge for $f$ between 0.009 and 0.01.
For $f>0.013$, the frequency response curve exhibits a shape resembling a linear system (although instabilities, not studied here, still persist).
%

The forcing amplitude range for which the IRC exists is the most troublesome from an engineering point of view, since at resonance the system response has two attractors with very different amplitudes.
The time series depicted in Fig.~\ref{numerical_results}b,c illustrate the phenomenon.
The two stable periodic solutions, separated by an unstable one, are shown in Fig.~\ref{numerical_results}d. In the figure, the shaded (clear) area is the basin of attraction of the smaller (larger) periodic solution.
We note that the basin of attraction of the IRC is very extended, which proves its practical relevance.
Very similar basins of attraction and overall dynamics were obtained for a single-DoF system with limited slip joint \cite{iwan1973transient}.

The direct numerical approach adopted so far is of course inadequate to analyze regions of existence of IRCs.
As discussed in the Introduction, several numerical and analytical techniques were developed for this aim.
In this study, singularity theory will be adopted, since it provides a well-established mathematical method to investigate the topology of curves in two dimensional spaces, as it is the case of frequency response curves.

\section{A brief review of singularity theory}

We adopt the framework developed in \cite{golubitsky1985} to study, through singularity theory, bifurcation diagrams defined by equations of the type \begin{equation}
g\left(x,\omega,\mu\right)=0,\label{g0}
\end{equation}
where $x$ is the state variable, $\omega$ is the bifurcation parameter (the frequency in the present study) and $\mu\in\mathbb{R}^n$ are additional parameters.
Here, we review only the basic elements of the theory necessary to understand its application to the identification of IRCs.
For a detailed exposition the reader is referred to \cite{golubitsky1985}.

Singularities are characterized by two groups of conditions on the derivatives of Eq.~(\ref{g0}) at a point.
Defining conditions, which set to zero the derivatives of (\ref{g0}), and nondegeneracy conditions, corresponding (in the simplest cases) to derivatives of (\ref{g0}) being different from zero.
The defining conditions always include \begin{equation}
g=\frac{\partial g}{\partial x}=0,\label{cond0}
\end{equation}
corresponding to the failure of the implicit function theorem.
Conditions in Eq.~(\ref{cond0}) alone indicate the presence of a fold, which is a persistent singularity, i.e. it is preserved by small perturbations.
To obtain a qualitative change in a bifurcation diagram, a nonpersistent singularity is necessary, i.e. a singularity that due to small perturbations (as, for example, a small change of a parameter value) disappears, leading to different possibilities.

Nonpersistent singularities require at least one additional defining condition.
The number of these additional defining conditions is the codimension of a singularity, which can be thought of as a measure of its complexity.
Codimension one singularities are the isola, the simple bifurcation and the hysteresis, whose defining and nondegeneracy conditions are \begin{itemize}
\item isola
\begin{equation}\label{eq:def_isola}
g=\frac{\partial g}{\partial \omega}=\frac{\partial g}{\partial x}=0,\quad\frac{\partial^2 g}{\partial x^2}\neq0,\quad\text{det}\left(\text d^2 g\right)>0,
\end{equation}
\item simple bifurcation 
\begin{equation}\label{eq:def_bif}
g=\frac{\partial g}{\partial \omega}=\frac{\partial g}{\partial x}=0,\quad\frac{\partial^2 g}{\partial x^2}\neq0,\quad\text{det}\left(\text d^2 g\right)<0,
\end{equation}
\item hysteresis 
\begin{equation}\label{eq:def_hys}
g=\frac{\partial g}{\partial x}=\frac{\partial^2 g}{\partial x^2}=0,\quad\frac{\partial g}{\partial \omega}\neq0,\quad\frac{\partial^3 g}{\partial x^3}\neq0,
\end{equation}
\end{itemize}
where $\text d^2g$ is the Hessian matrix of $g\left(x,\omega\right)$.

\begin{figure}[tb]
\begin{centering}
\setlength{\unitlength}{\textwidth}
\begin{picture}(1,0.95)
\put(0.0,0.29){\includegraphics[trim = 15mm 10mm 12mm 10mm,clip,width=0.3\textwidth]{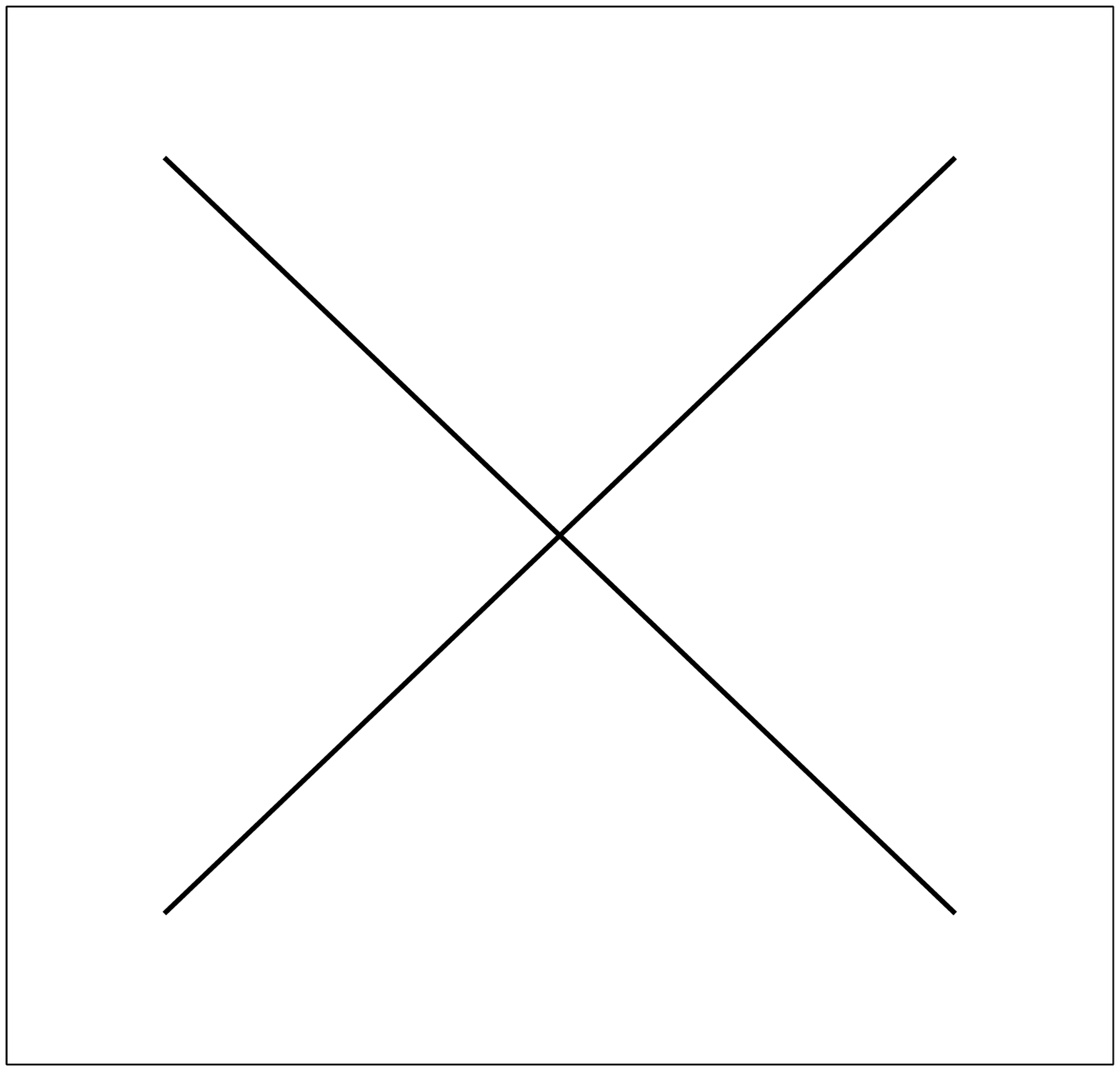}}
\put(0.36,.29){\includegraphics[trim = 15mm 10mm 12mm 10mm,clip,width=0.3\textwidth]{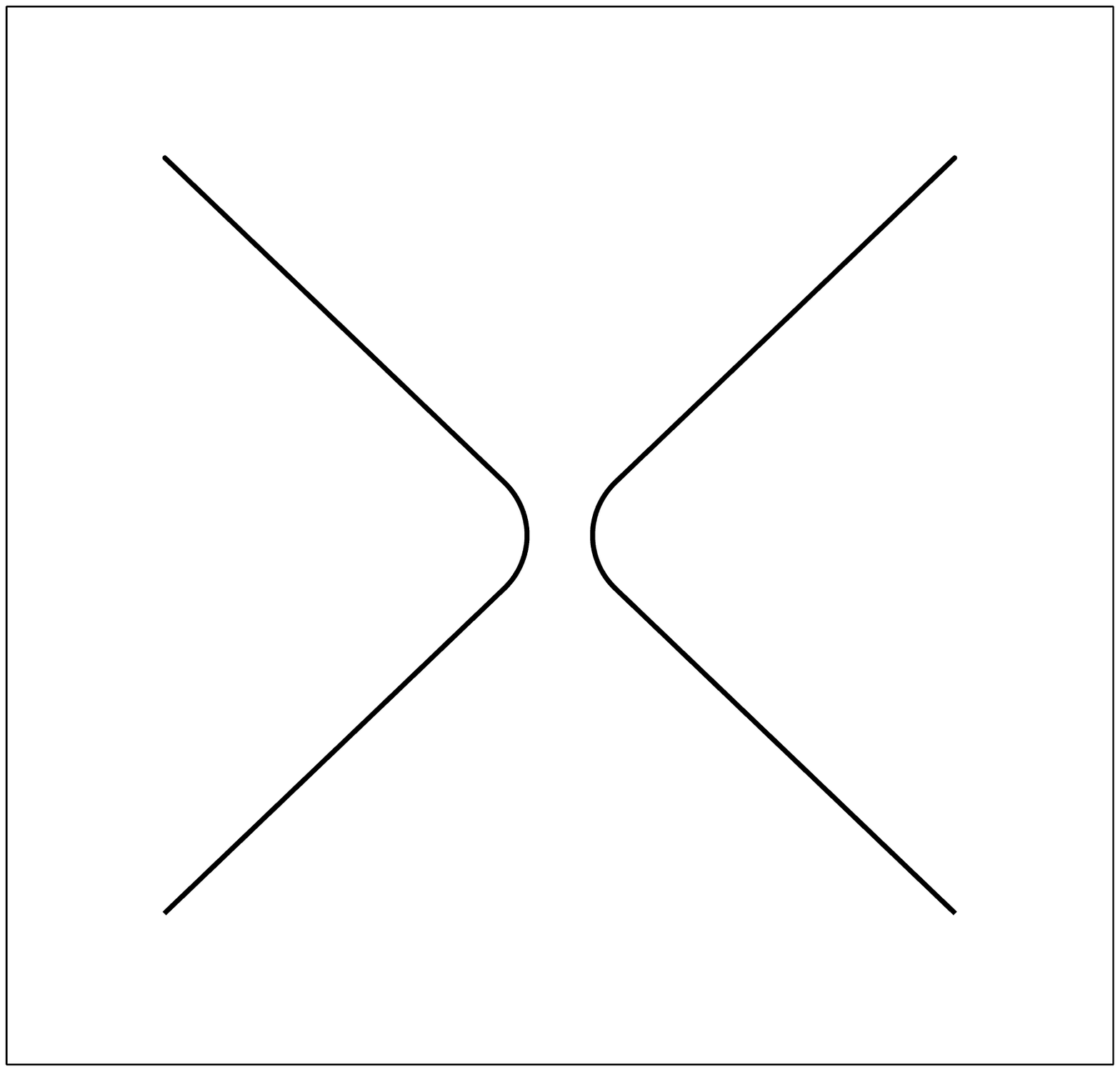}}
\put(0.67,0.29){\includegraphics[trim = 15mm 10mm 12mm 10mm,clip,width=0.3\textwidth]{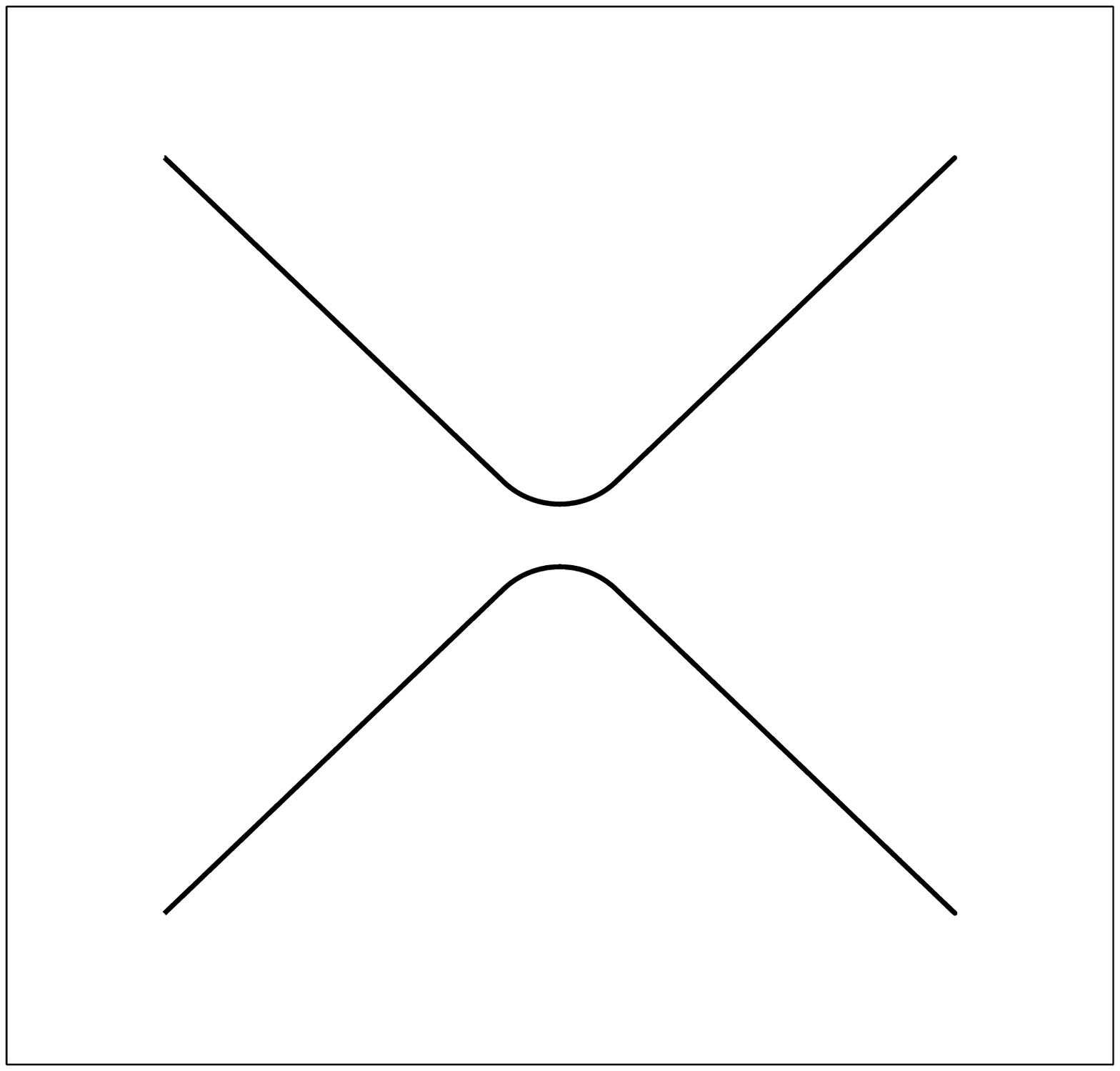}}
\put(0.0,-0.01){\includegraphics[trim = 15mm 10mm 12mm 10mm,clip,width=0.3\textwidth]{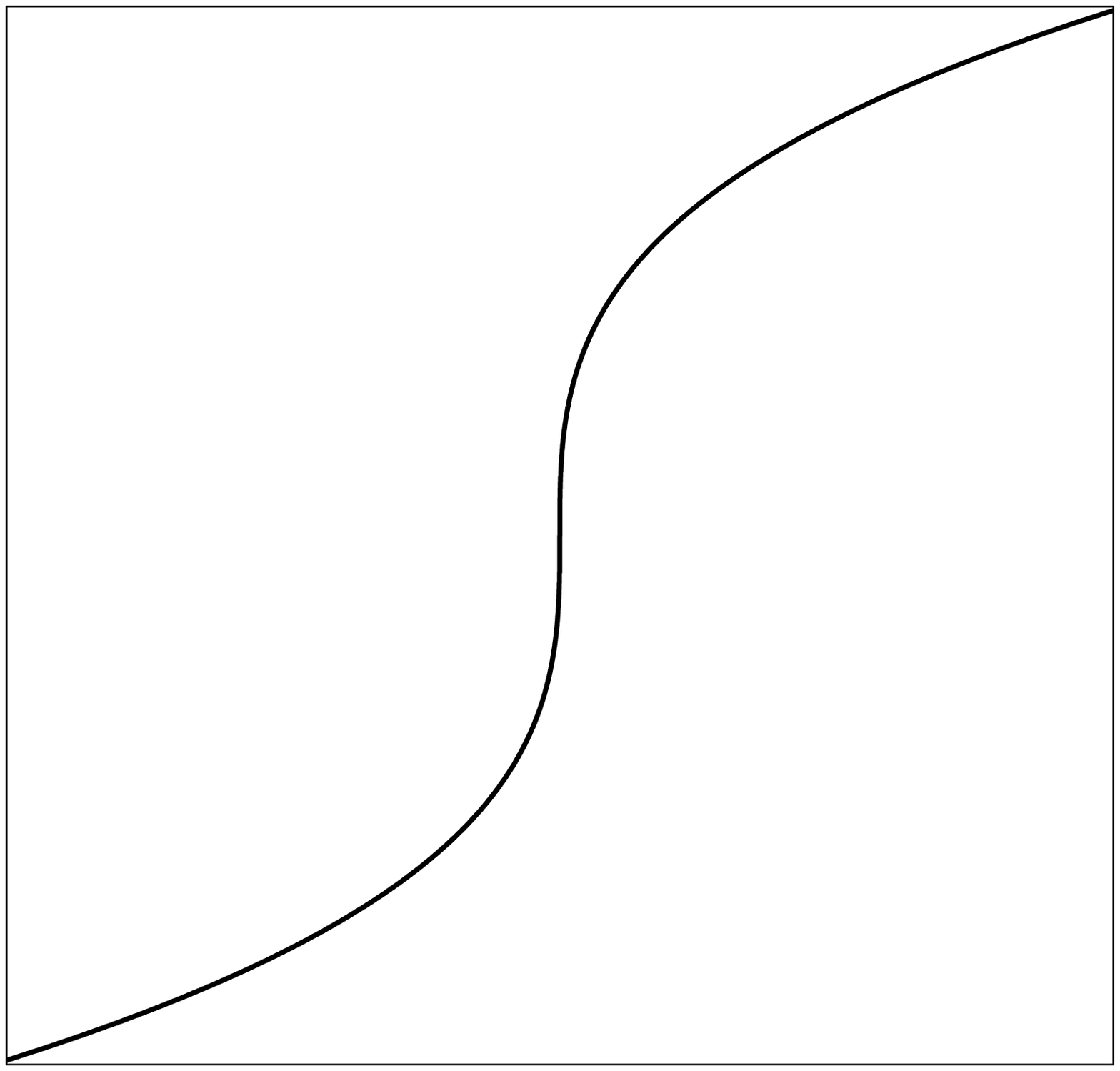}}
\put(0.36,-0.01){\includegraphics[trim = 15mm 10mm 12mm 10mm,clip,width=0.3\textwidth]{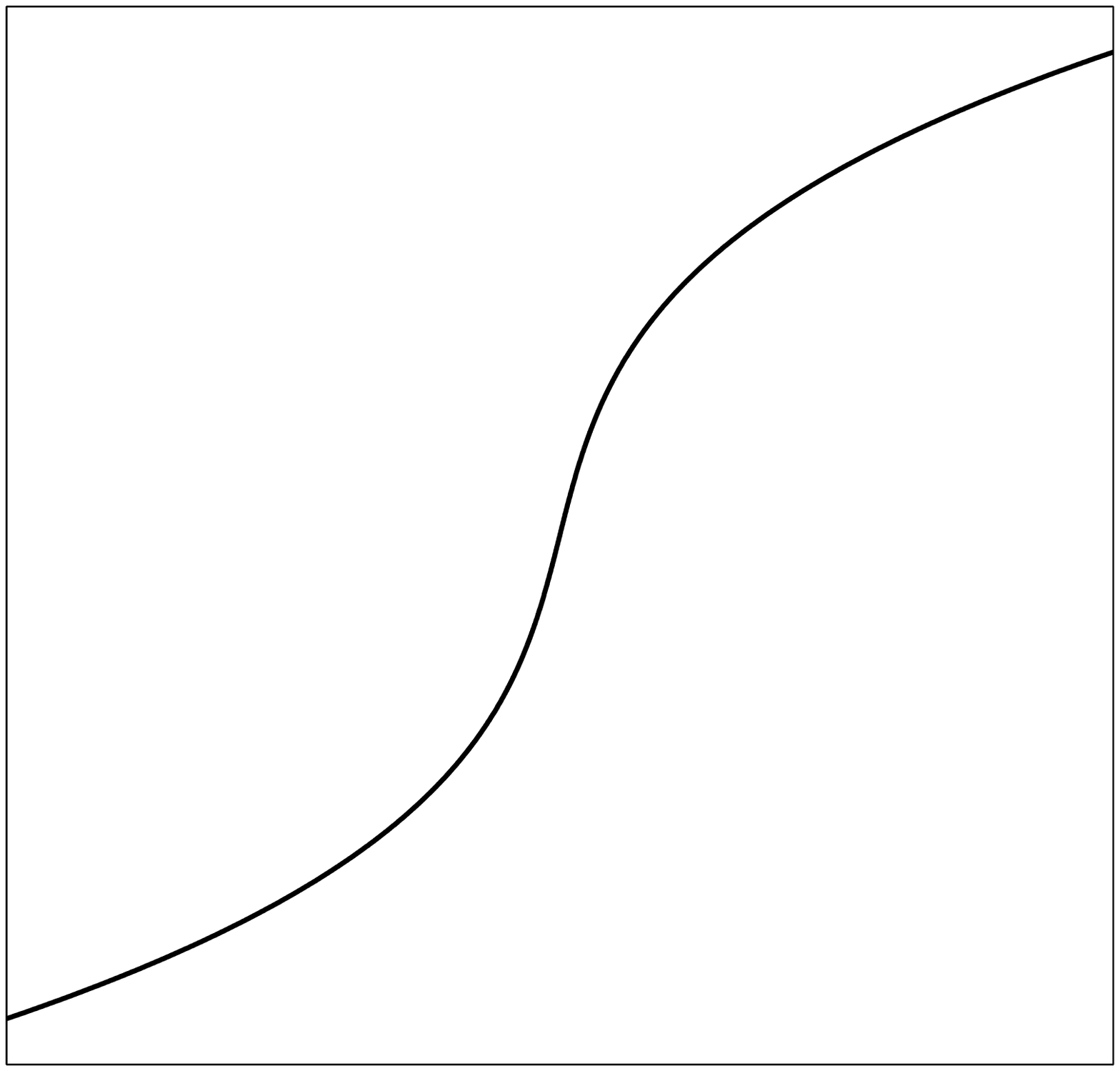}}
\put(0.67,-.01){\includegraphics[trim = 15mm 10mm 12mm 10mm,clip,width=0.3\textwidth]{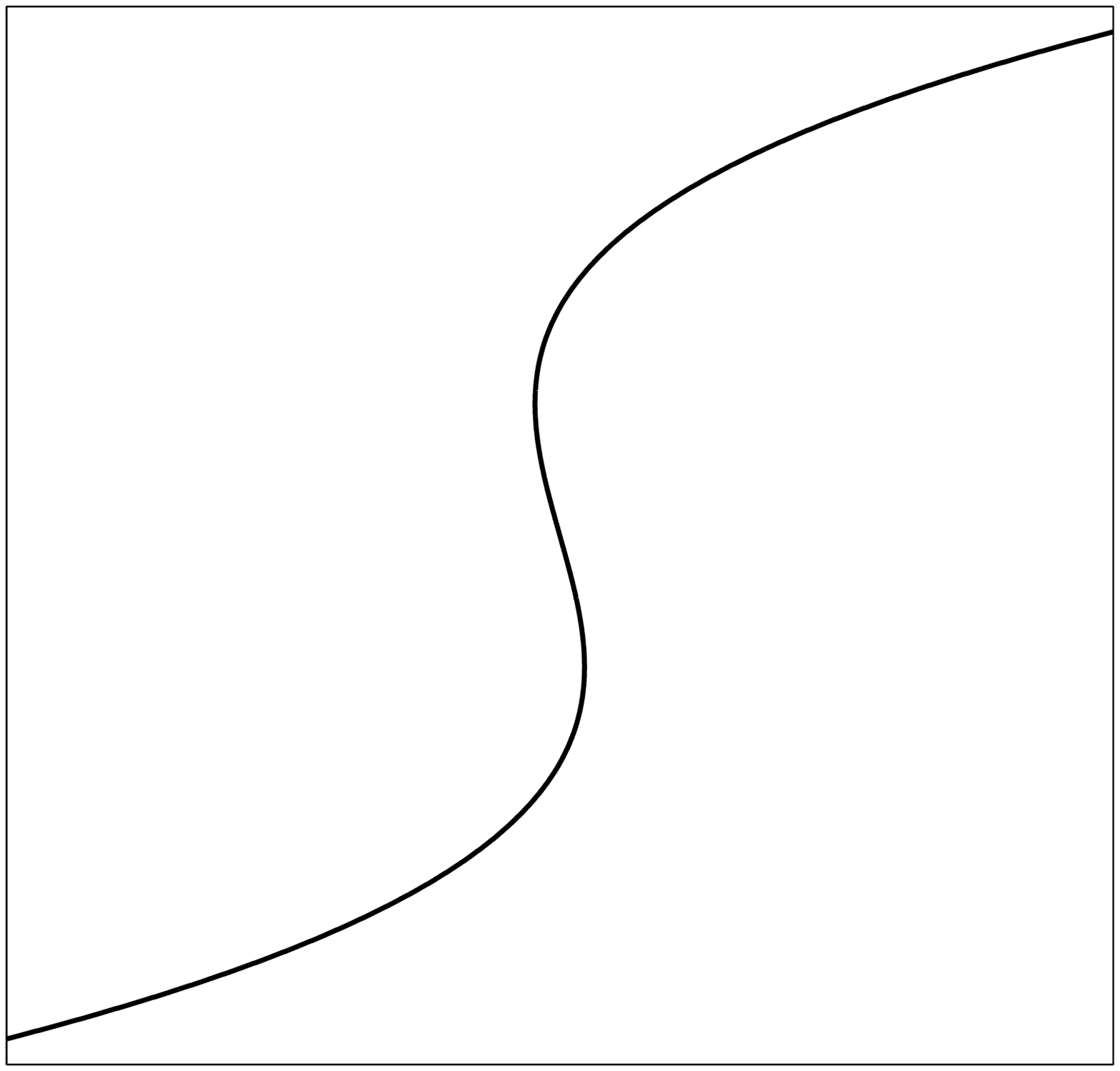}}
\put(0.0,.59){\includegraphics[trim = 15mm 10mm 12mm 10mm,clip,width=0.3\textwidth]{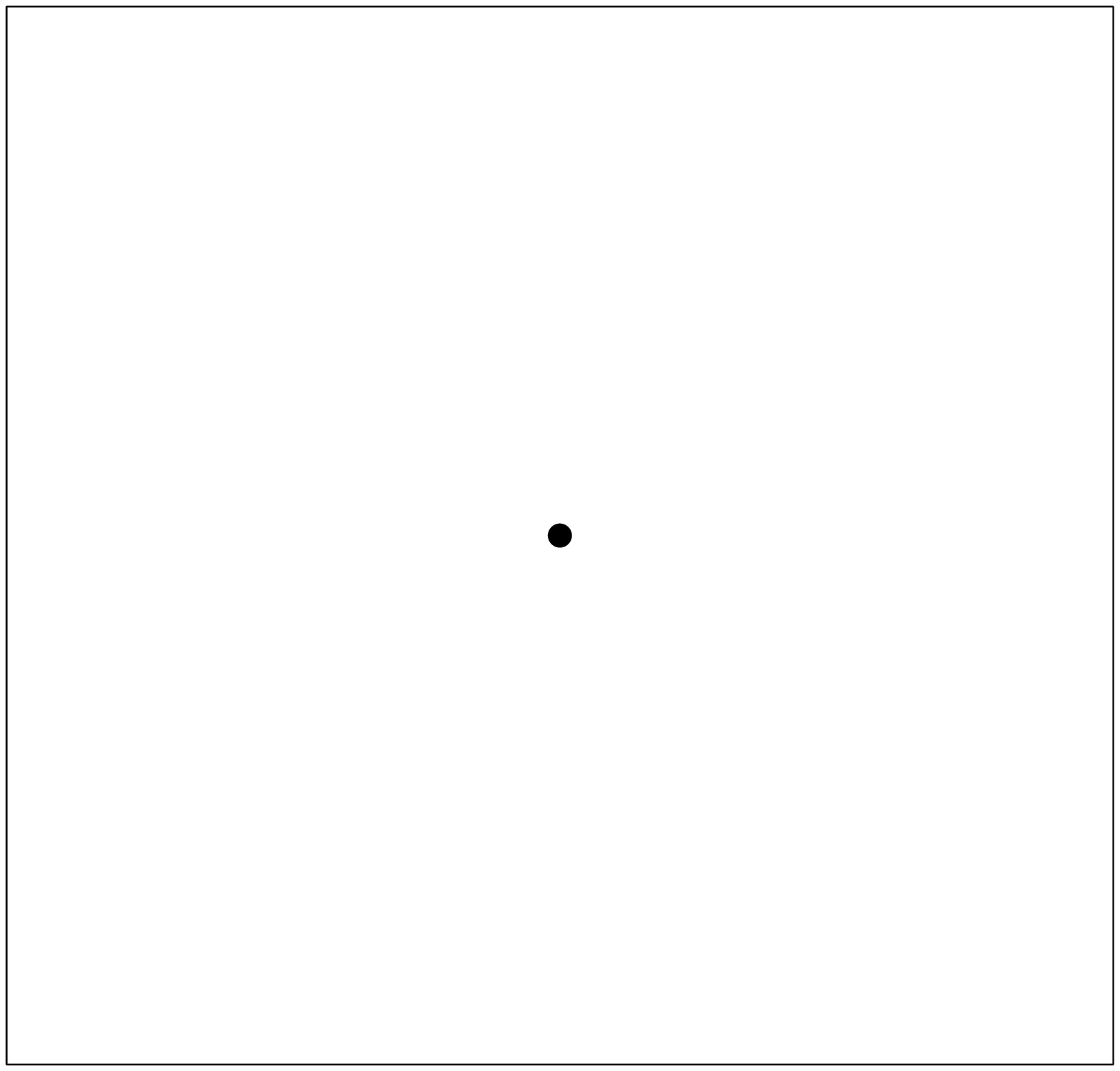}}
\put(0.36,.59){\includegraphics[trim = 15mm 10mm 12mm 10mm,clip,width=0.3\textwidth]{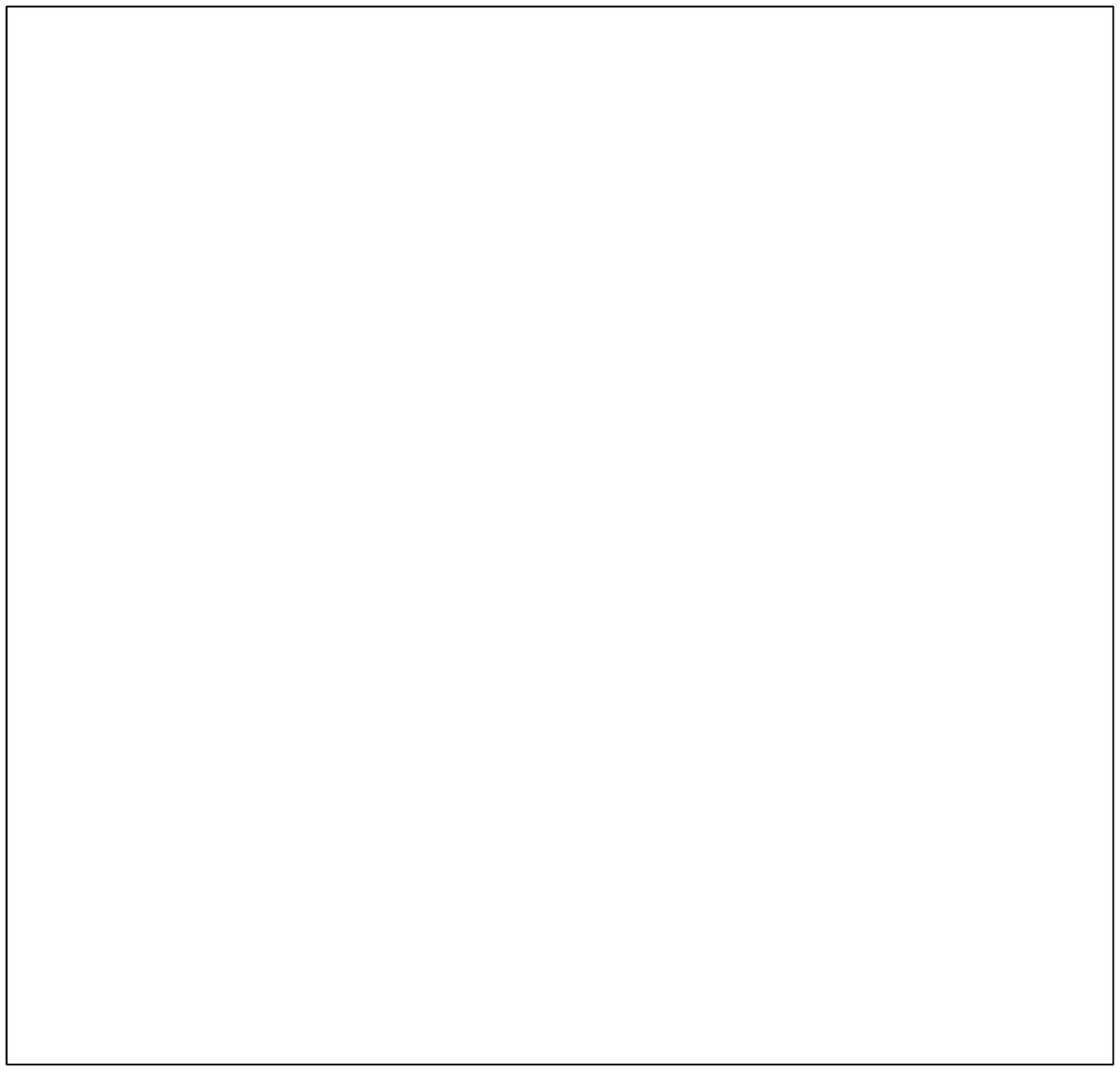}}
\put(0.67,.59){\includegraphics[trim = 15mm 10mm 12mm 10mm,clip,width=0.3\textwidth]{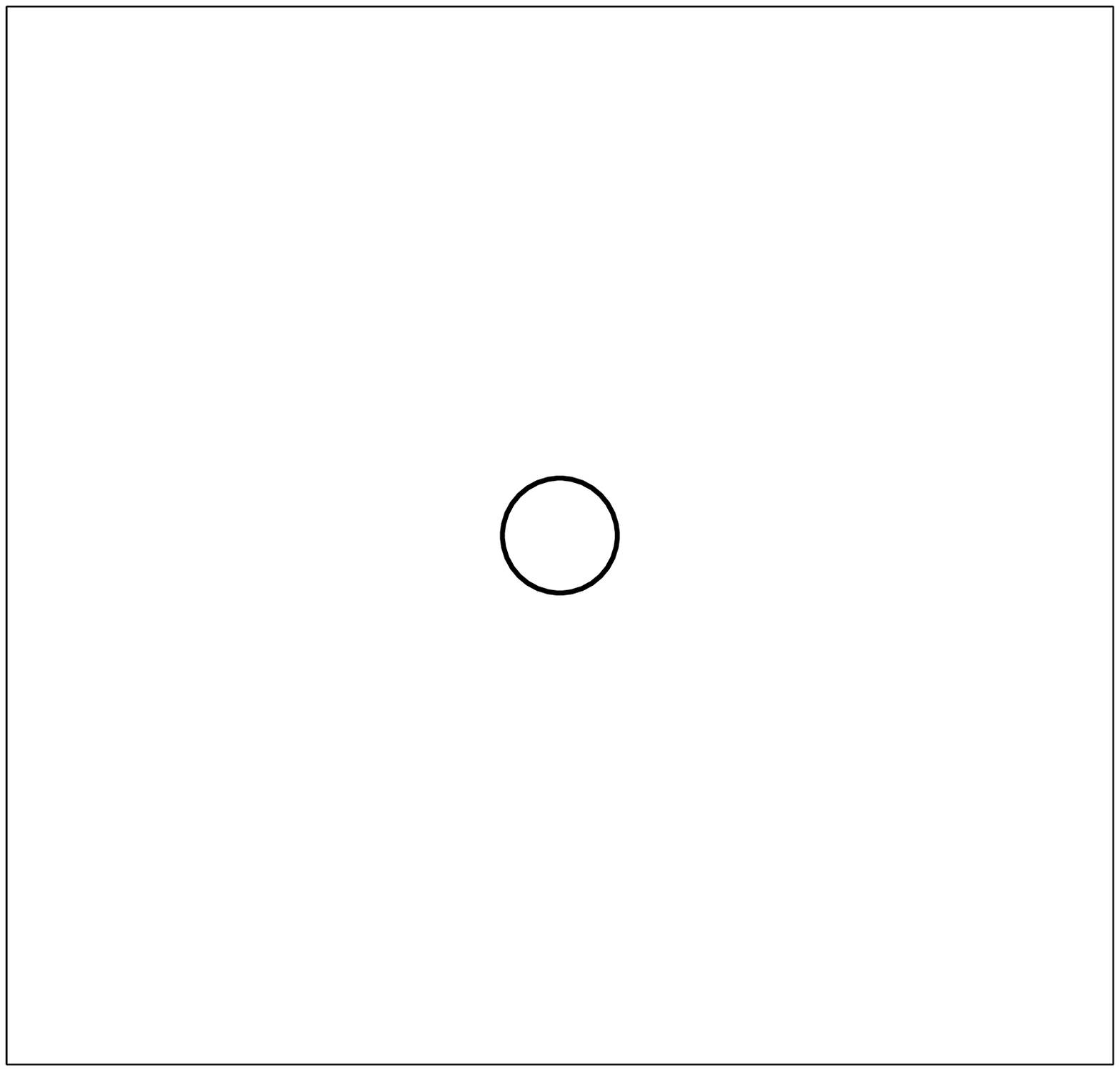}}
\put(0.03,0.835){\textbf{(a)}}
\put(0.03,0.535){\textbf{(d)}}
\put(0.03,0.235){\textbf{(g)}}
\put(0.39,0.835){\textbf{(b)}}
\put(0.39,0.535){\textbf{(e)}}
\put(0.39,0.235){\textbf{(h)}}
\put(0.7,0.835){\textbf{(c)}}
\put(0.7,0.535){\textbf{(f)}}
\put(0.7,0.235){\textbf{(i)}}
\put(0.045,0.89){Unperturbed diagrams}
\put(0.56,0.89){Perturbed diagrams}
\put(0.335,0.01){\linethickness{.3mm}\line(0,1){0.9}}
\end{picture}
\end{centering}
\caption{(a-c) The isola and its perturbations; (d-f) the simple bifurcation and its perturbations; (g-i) the hysteresis and its perturbations.}\label{fig:singularities}
\end{figure}

Isola singularities are found when an IRC appears. This corresponds to the presence of an isolated solution as in Fig.~\ref{fig:singularities}a. When perturbed, this nonpersistent diagram can only result in two outcomes: either no solutions (Fig.~\ref{fig:singularities}b) or a closed branch of solutions (Fig.~\ref{fig:singularities}c).

Simple bifurcation points correspond to the center of X-shaped diagrams, as illustrated in Fig.~\ref{fig:singularities}d.
The corresponding perturbations are shown in  Fig.~\ref{fig:singularities}e,f.
This is the local phenomenon that underlies the merging of an IRC with a main branch.

Finally, an example of hysteresis is depicted in Fig.~\ref{fig:singularities}g.
Geometrically, it is an inflexion point characterized by a vertical tangent.
When perturbed, a hysteresis results in one of the diagrams illustrated in Fig.~\ref{fig:singularities}h,i.

As it can be intuitively understood, nonpersistent singularities act as boundaries between qualitatively different persistent diagrams.
Persistence is related to the concept of equivalence: a diagram is persistent if sufficiently small perturbations result in equivalent diagrams, otherwise it is nonpersistent.
For a rigorous definition of equivalence we refer the interested reader to \cite{golubitsky1985}.
Thanks to this property, singularities can be used as a tool to divide the parameter space in different zones, characterized by a unique (up to equivalence) response.

Throughout the present study, we take advantage of several nonpersistent singularities to identify appearance and merging of IRCs. Conditions for the existence of isola and simple bifurcation singularities are used.
This will allow us to define regions in parameter space where IRCs exist.

\section{Analytical development}

To obtain an algebraic equation characterizing the frequency response of system (\ref{NL_DF}), we adopt a harmonic balance procedure, approximating steady state solutions by 
\begin{equation}
  x=\gamma e^{j\omega t}+\bar \gamma e^{-j\omega t},
\end{equation}
where $j$ is the imaginary unit and the bar indicates complex conjugates.
Neglecting higher-order harmonics, we have 
\begin{equation}
  -\gamma  \omega ^2+\gamma +j \gamma  c_1 \omega +3 j \gamma ^2 c_3 \omega ^3 \bar \gamma +10 j \gamma ^3 \omega ^5 \left(\bar \gamma\right)^2-f=0.
\end{equation}
We introduce polar coordinates $\gamma=ae^{j\phi}$, with $a$ and $\phi$ real, obtaining
\begin{equation}\label{eq:imagreal}
  10 j a^5 \omega ^5+3 j a^3 c_3 \omega ^3+j a c_1 \omega -a \omega ^2+a=f e^{-j\phi },
\end{equation}
separating then real and imaginary parts of Eq.~(\ref{eq:imagreal}) and summing up their square, yields the equation
\begin{equation}\label{eq:g_damp}
g\left(A,\Omega,F,c_1,c_3\right)=100 A^5 \Omega ^5+60 A^4 c_3 \Omega ^4+A^3 \Omega ^3 \left(20 c_1+9 c_3^2\right)+6 A^2 c_1 c_3 \Omega ^2+A \left(\left(c_1^2-2\right) \Omega +\Omega ^2+1\right)-F=0,
\end{equation}
where $A=a^2$, $\Omega=\omega^2$ and $F=f^2$.
Solutions of~\eqref{eq:g_damp} approximate the frequency response of the system under study.

To obtain expressions corresponding to appearance and merging of IRCs, we use the conditions characterizing the corresponding singularities, isola and simple bifurcation (Eqs.~(\ref{eq:def_isola}) and (\ref{eq:def_bif})).

Starting from Eq.~\eqref{eq:g_damp} we obtain
\begin{equation}\label{eq:derOmega1}
\frac{\partial g}{\partial A}=500 A^4 \Omega ^5+240 A^3 c_3 \Omega ^4+3 A^2 \Omega ^3 \left(20 c_1+9 c_3^2\right)+\Omega ^2 (12 A c_1 c_3+1)+\left(c_1^2-2\right) \Omega +1
\end{equation}
and 
\begin{equation}\label{eq:derA1}
  \frac{\partial g}{\partial \Omega}=A \left(500 A^4 \Omega ^4+240 A^3 c_3 \Omega ^3+27 A^2 c_3^2 \Omega ^2+12 A c_1 \Omega  (5 A \Omega +c_3)+c_1^2+2 \Omega -2\right).
\end{equation}
While
\begin{equation}
  \frac{\partial^2 g}{\partial A^2}=2 \Omega ^2 \left(A \Omega  \left(1000 A^2 \Omega ^2+360 A c_3 \Omega +27 c_3^2\right)+6 c_1 (10 A \Omega +c_3)\right)
\end{equation}
and 
\begin{equation}
  \begin{split}
    \text{det}\left(\text d^2 g\right)=&-2250000 A^8 \Omega ^8-1920000 A^7 c_3 \Omega ^7-592200 A^6 c_3^2 \Omega ^6-240 A^4 \Omega ^5 \left(324 A c_3^3+25\right)\\
    &-5 A^3 \Omega ^4 \left(729 A c_3^4-2000 A+480 c_3\right)+4 \Omega ^2 \left(81 A^2 c_3^2-1\right)+24 A^2 c_3 \Omega ^3 (160 A-9 c_3)\\
    &-2 c_1^2 \left(11500 A^4 \Omega ^4+3840 A^3 c_3 \Omega ^3+297 A^2 c_3^2 \Omega ^2+2 \Omega -2\right)\\
    &-24 A c_1 \Omega  \left(c_3 \left(10200 A^4 \Omega ^4+3 \Omega -4\right)+10 A \Omega  \left(1750 A^4 \Omega ^4+2 \Omega -3\right)+1875 A^3 c_3^2 \Omega ^3+108 A^2 c_3^3 \Omega ^2\right)\\
    &-24 A c_1^3 \Omega  (15 A \Omega +2 c_3)-c_1^4+8 \Omega -4.
  \end{split}
\end{equation}
$g=\partial g/\partial \Omega=\partial g/\partial A=0$ forms a nonlinear algebraic system of equations with the unknown $F$, $A$ and $\Omega$, whose real solutions are
\begin{align}
S_1:\,\left(F_1,A_1,\Omega_1\right)&=\left(0,\frac{1}{20} \left(-3 c_3-\sqrt{9 c_3^2-40 c_1}\right),1\right)\\
S_2:\,\left(F_2,A_2,\Omega_2\right)&=\left(0,\frac{1}{20} \left(-3 c_3+\sqrt{9 c_3^2-40 c_1}\right),1\right)\\
S_3:\,\left(F_3,A_3,\Omega_3\right)&=\left(F_3,\frac{1}{100} \left(-9 c_3-\sqrt{81 c_3^2-200 c_1}\right),1\right)\\
S_4:\,\left(F_4,A_4,\Omega_4\right)&=\left(F_4,\frac{1}{100} \left(-9 c_3+\sqrt{81 c_3^2-200 c_1}\right),1\right),
\end{align}
where
\begin{equation}
  \begin{split}
    F_3&=\frac{\left(-9 c_3-\sqrt{81 c_3^2-200 c_1}\right) \left(200 c_1-3 c_3 \left(\sqrt{81 c_3^2-200 c_1}+9 c_3\right)\right)^2}{6250000}\\
    F_4&=\frac{\left(-9 c_3+\sqrt{81 c_3^2-200 c_1}\right) \left(200 c_1+3 c_3 \left(\sqrt{81 c_3^2-200 c_1}-9 c_3\right)\right)^2}{6250000}.
  \end{split}
\end{equation}
Remarkably, all of them appear at $\Omega=1$, clearly due to the linear nature of the elastic force.
Solutions are physically meaningful only if $A$, $\Omega$ and $F$ are real and non-negative.
We call $S_1$ the singularity relative to ($A_1,\Omega_1,F_1$) and so on for the others.

Recalling the assumption that linear damping is positive ($c_1>0$) we have
\begin{equation}\label{eq:A1pos}
  A_1\geq0\,\&\,A_2\geq0\quad\iff\quad c_3\leq c_{31}=-\frac{2}{3}\sqrt{10c_1}
\end{equation}
while
\begin{equation}\label{eq:A3pos}
A_3\geq0\,\&\,A_4\geq0\,\&\,F_3\geq0\,\&\,F_4\geq0\quad\iff\quad c_3\leq c_{32}=-\frac{10}{9}\sqrt{2c_1}
\end{equation}
(it can be verified that, for $\Omega=1$, $F>0$ if and only if $A>0$, unless $A=A_1$ or $A=A_2$).
From these results and the fact that $c_{32}>c_{31}$, it follows that, if $c_3>c_{32}$, Eq.~\eqref{eq:g_damp} does not predict the existence of any singularity of the studied type.

We call $H_i=\text{det}\left(\text d^2 g\right)|_{A=A_i,\Omega=1}$, $i=1...4$.
These quantities verify
\begin{itemize}
\item $H_1>0\,\iff\, c_3<c_{31}$; $H_1=0\,\iff\,c_3=c_{31}$
\item $H_2>0\,\iff\, c_3<c_{31}$; $H_2=0\,\iff\,c_3=c_{31}$
\item $H_1>0\,\iff\, c_3<c_{32}$; $H_3=0\,\iff\,c_3=c_{32}$
\item $H_4>0\,\iff\, c_{31}<c_3<c_{32}$;  $H_4<0\,\iff\, c_3<c_{31}$; $H_4=0\,\iff\,c_3=c_{31}$ or $c_3=c_{32}$.
\end{itemize}
The case of $c_3>c_{32}$ and conditions for which $H_i$, $i=1...4$ are not considered.

\begin{figure}
\begin{centering}
\setlength{\unitlength}{\textwidth}
\begin{picture}(1,0.5)
\put(0.03,-0.01){\includegraphics[trim = 10mm 10mm 12mm 10mm,clip,width=0.45\textwidth]{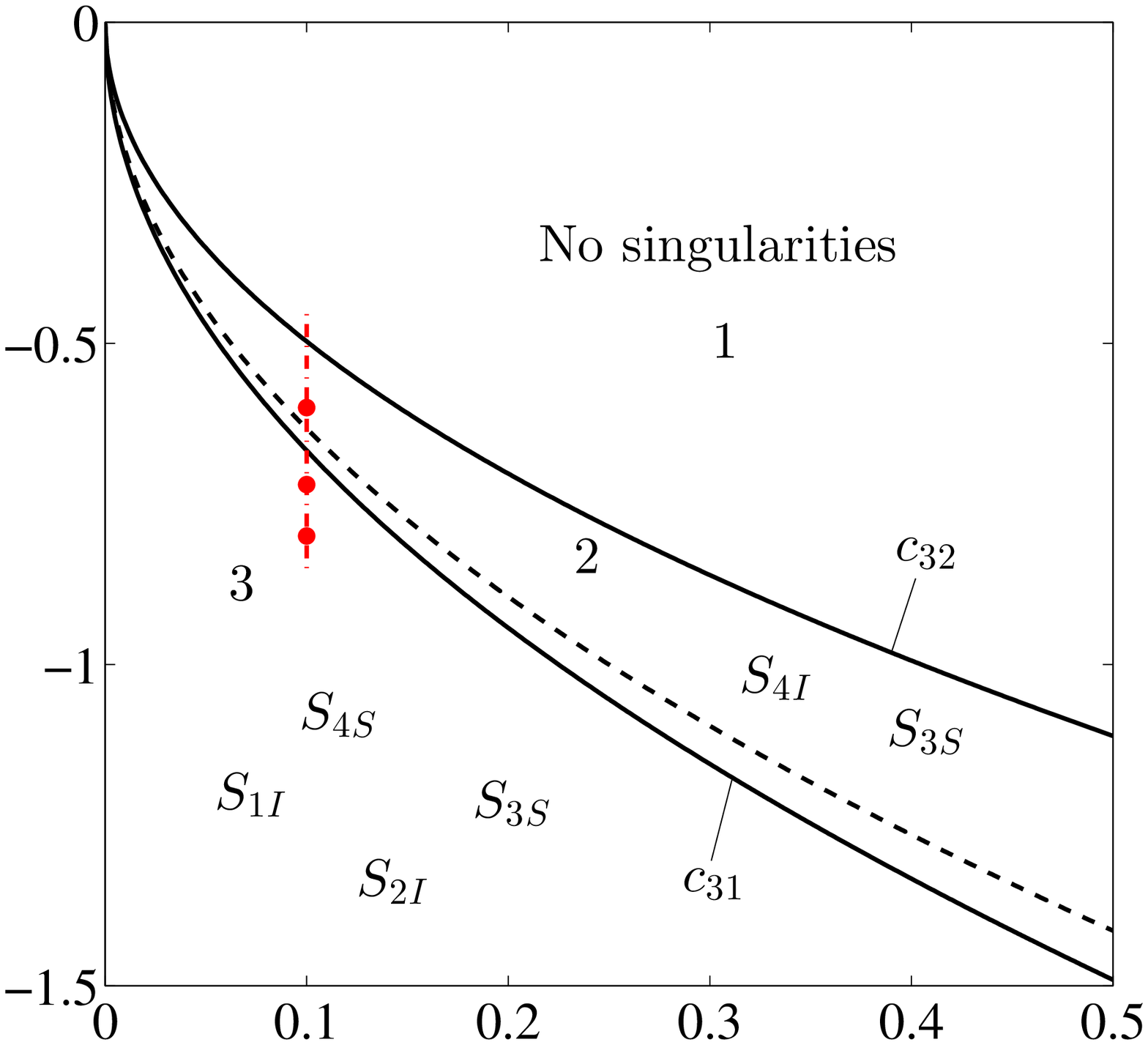}}
\put(0.52,-0.01){\includegraphics[trim = 10mm 10mm 12mm 10mm,clip,width=0.45\textwidth]{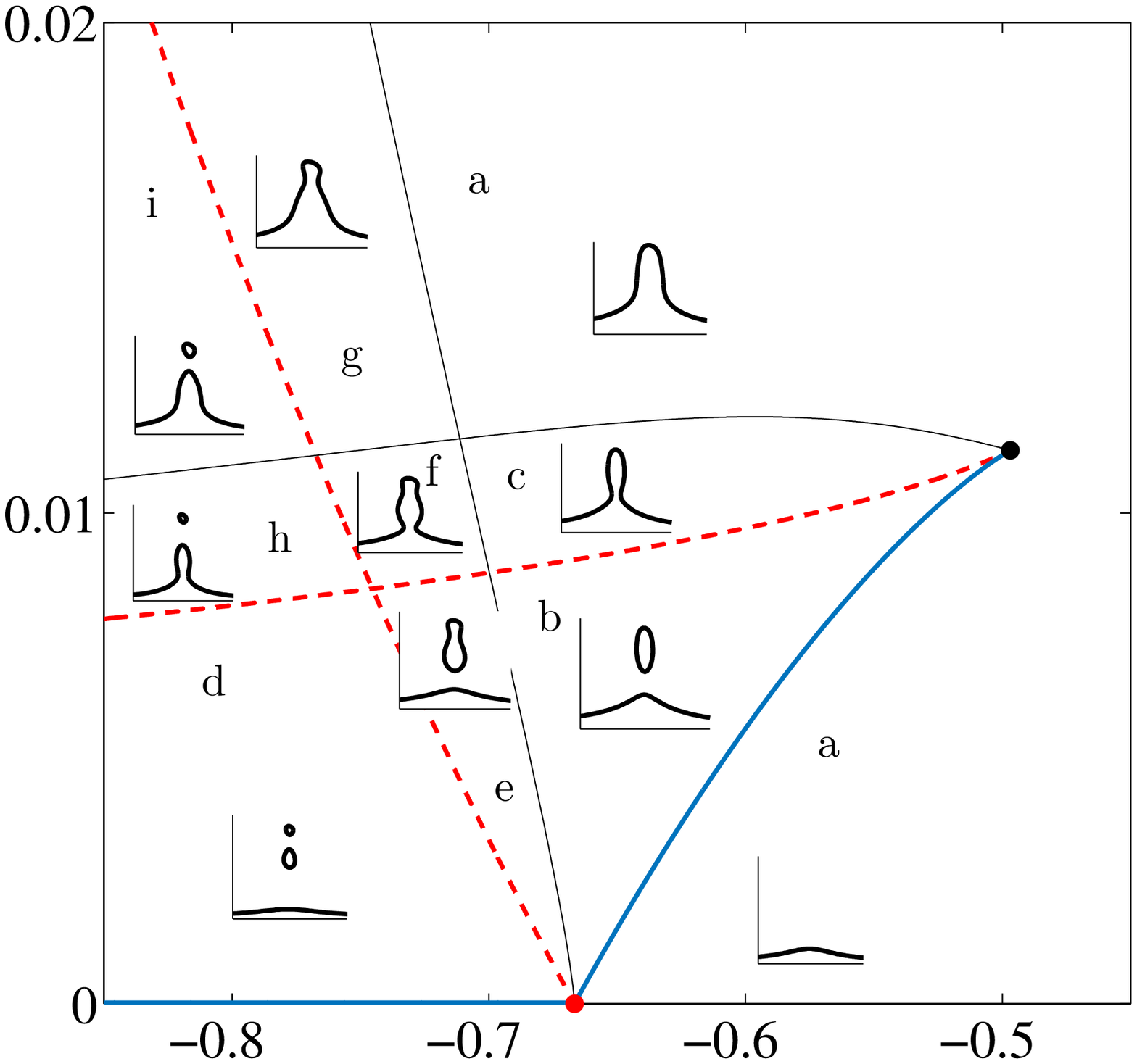}}
\put(0.075,0.37){\textbf{(a)}}
\put(0.565,0.37){\textbf{(b)}}
\put(0.255,-0.03){$c_1$}
\put(0.75,-0.03){$c_3$}
\put(0,0.2){\rotatebox{90}{$c_3$}}
\put(0.49,0.2){\rotatebox{90}{$f$}}
\end{picture}
\end{centering}
\caption{(a) Region of existence of the studied singularities in the $c_1,c_3$ space. Subscript S indicates ``simple bifurcation", while subscript I indicates ``isola". (b) Regions exhibiting different frequency response scenarios in the $c_3,f$ parameter space for $c_1=0.1$. Solid blue lines: isola singularities, dashed red lines: simple bifurcations, thin black lines: hysteresis singularities; black and red dots mark higher codimension singularities.\label{region_c1c3}}
\end{figure}

Thus, taking into account Eqs.~\eqref{eq:def_isola} and \eqref{eq:def_bif} and the conditions in Eqs.~\eqref{eq:A1pos} and \eqref{eq:A3pos}, we can conclude the following:
\begin{itemize}
\item $S_1$ and $S_2$ exist and they are isolas $\iff\, c_3<c_{31}$
\item $S_3$ exists and it is a simple bifurcation $\iff\,c_3<c_{32}$
\item $S_4$ exists and it is an isola $\iff\,c_{31}<c_3<c_{32}$
\item $S_4$ exists and it is a simple bifurcation $\iff\,c_3<c_{31}$
\end{itemize}
Fig.~\ref{region_c1c3}a summarizes the scenario described above.

Note that $S_1$ and $S_2$ correspond to isola generated at $F=0$.
Such motion is possible, since for $c_3<-2\sqrt{c_1}$ the damping force becomes negative in the interval $\sqrt{\left(-c_3-\sqrt{c_3^2-4 c_1}\right)/2}<\left|\dot x\right|<\sqrt{\left(-c_3+\sqrt{c_3^2-4 c_1}\right)/2}$,
and the system can develop periodic motions even in absence of external forcing.
The dashed line in Fig.~\ref{region_c1c3} marks the region having negative damping force.
It can be noticed that $c_{31}$ lays below this line, therefore $S_1$ and $S_2$ exist only if the damping force is negative for a certain velocity range.

For $c_3=c_{32}$, $S_4$ (isola singularity) and $S_3$ (simple bifurcation) merge in a more degenerate singularity.
At this point, besides the conditions defining isolas and simple bifurcations, two other conditions are verified so that
\begin{equation}
  \label{eq:def-wc}
  \frac{\partial g}{\partial A}=\frac{\partial^2 g}{\partial A^2}=\frac{\partial g}{\partial \Omega}=\frac{\partial^2 g}{\partial A\partial\Omega}=0.
\end{equation}
These, joined with the nondegeneracy conditions
\begin{equation}
  \label{eq:non-cusp}
  \left.\frac{\partial^3 g}{\partial A^3}\right|_{A=A_3=A_4, \Omega=1, c_3=c_{32}}=\frac{160}{3}c_1\neq 0, \qquad \left.\frac{\partial^2 g}{\partial \Omega^2}\right|_{A=A_3=A_4, \Omega=1, c_3=c_{32}}=\frac{\sqrt{2 c_1}}{5}\neq0,
\end{equation}
correspond to a codimension three singularity, the winged cusp.

Line $c_3=c_{31}$ is the locus of an even more degenerate singularity.
For $c_3=c_{31}$, we have $A_1=A_2=A_4$, $F_1=F_2=F_4=0$, and 
\begin{equation}
  \label{eq:singge4}
  \frac{\partial g}{\partial A}=\frac{\partial^2 g}{\partial A^2}=\frac{\partial^3 g}{\partial A^3}=\frac{\partial g}{\partial \Omega}=\frac{\partial^2 g}{\partial A\partial \Omega}=0, \quad \frac{\partial^4 g}{\partial A^2}\neq 0,
\end{equation}
which corresponds to a singularity of codimension $\geq4$.
Considering its complexity and the fact that several additional parameters should be introduced to describe all the existing persistent diagrams in its neighborhood, we let its detailed analysis to future studies. We only mention that, since this singularity appears at $F=0$, it is possible that part of its complexity is due to phenomena appearing in the case $F<0$ that does not have a physical meaning.

Zones 1, 2 and 3 of Fig.~\ref{region_c1c3}a correspond to qualitatively different behaviors of the system.
In zone 1 the system has no singularities, therefore a behavior topologically similar to that of a linear system can be expected. In zone 2 an isola singularity and a simple bifurcation are present. Correspondingly, an IRC appears at a certain forcing amplitude, which then merges with the main branch for a sufficiently high value of $f$.
In zone 3, where two isola singularities and two simple bifurcation singularities exist, an even more involved scenario is expected.

To obtain further insight into the dynamics of the system, we fix the parameter $c_1$ at 0.1 and plot the singularities of the system in the $c_3,f$ parameter space. The result is illustrated in Fig.~\ref{region_c1c3}b, where hysteresis singularities are also represented.
Solid blue lines and dashed red lines indicate isola and simple bifurcation singularities, respectively.
The black dot indicates a winged cusp while the red dot marks the undefined higher codimension singularity.
Zone 1 of Fig.~\ref{region_c1c3}a is on the right of the winged cusp, for $c_3>-0.498$. No relevant phenomena is noticed here, because of the absence of singularities.
For $-0.667<c_3<-0.498$, between the red and the black dots, the system is in zone 2. For increasing values of $f$, at the passage between regions a and b of Fig.~\ref{region_c1c3}b, an IRC appears, which then merges with the main branch as the system enters in region c. Successively, passing the line of hysteresis singularities and reaching again region a, the frequency response becomes again monovalued (excluding non-periodic solutions, not analyzed in this study).
This transition corresponds to the one described in Fig.~\ref{numerical_results}a for $c_1=0.1$ and $c_3=-0.6$.
According to the analytical computation, appearance of the IRC and its merging should occur at $f=0.0056$ and $f=0.0096$ at an amplitude of $x=0.58$ and $x=0.30$, respectively.
The prediction is confirmed numerically, as illustrated in Fig.~\ref{numerical_results}a.

\begin{figure}
\begin{centering}
\setlength{\unitlength}{\textwidth}
\begin{picture}(1,0.5)
\put(0.03,-0.01){\includegraphics[trim = 10mm 10mm 12mm 10mm,clip,width=0.45\textwidth]{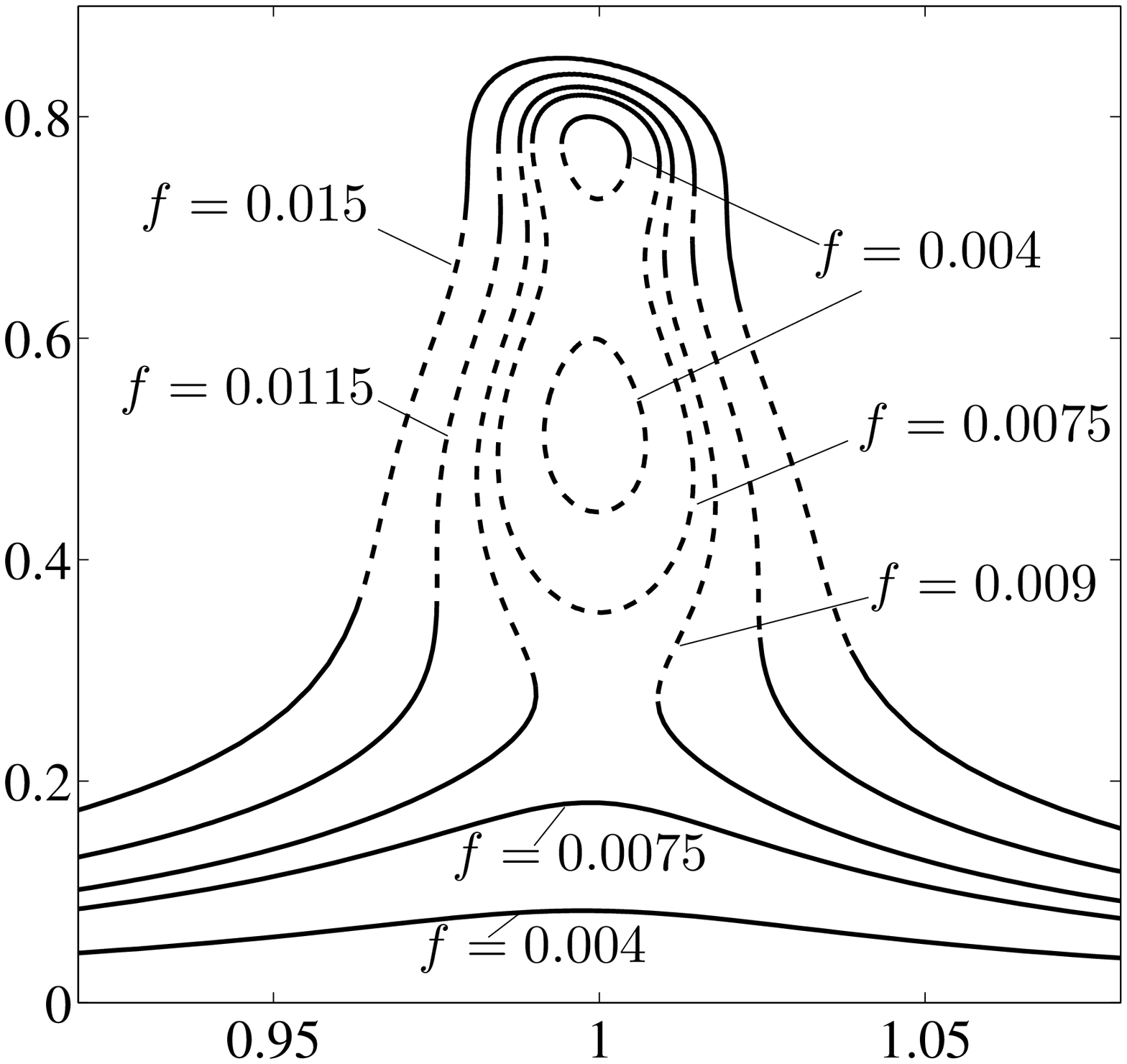}}
\put(0.52,-0.01){\includegraphics[trim = 10mm 10mm 12mm 10mm,clip,width=0.45\textwidth]{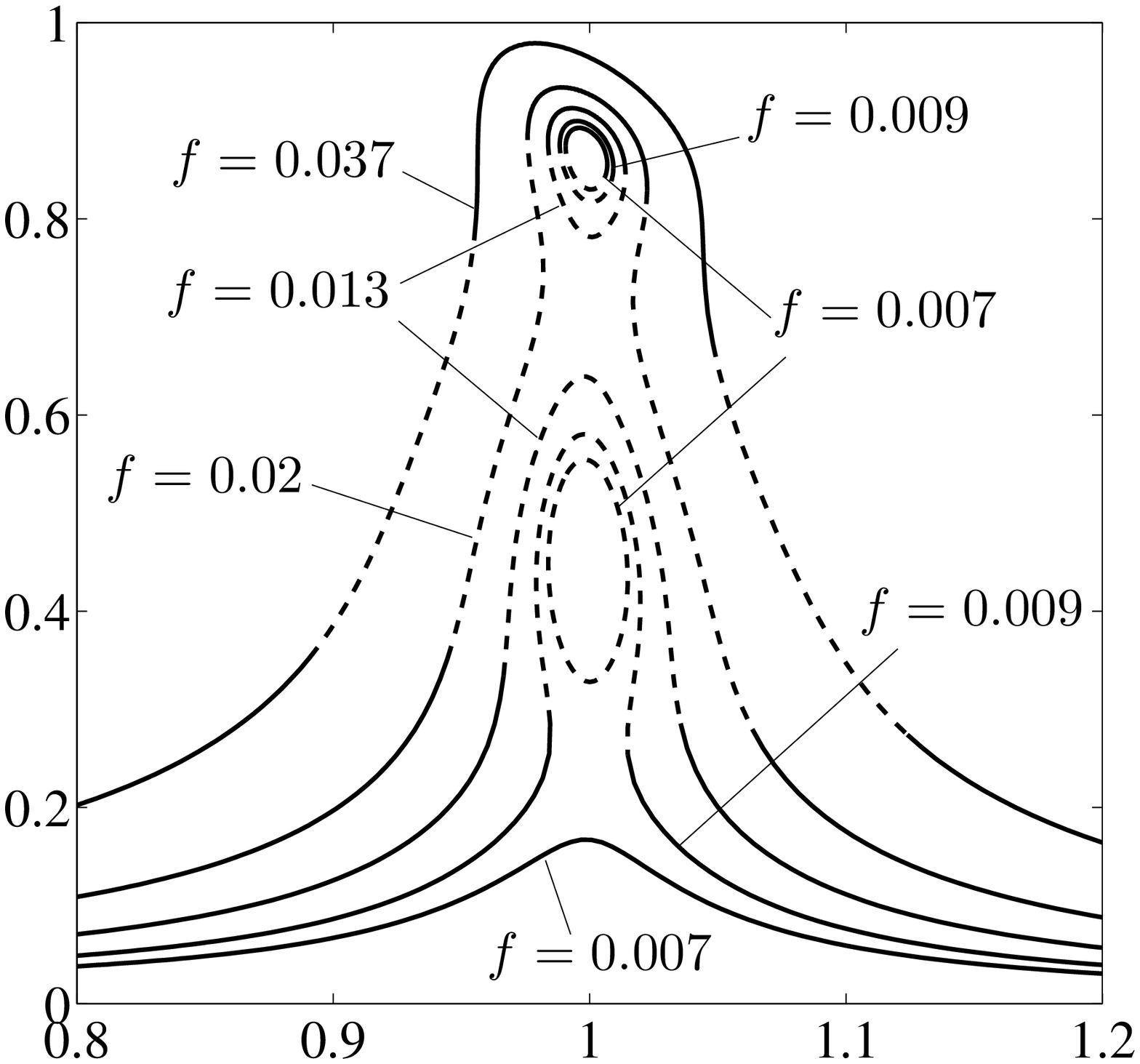}}
\put(0.075,0.37){\textbf{(a)}}
\put(0.565,0.37){\textbf{(b)}}
\put(0.255,-0.03){$\omega$}
\put(0.75,-0.03){$\omega$}
\put(0,0.2){\rotatebox{90}{$x$}}
\put(0.49,0.2){\rotatebox{90}{$x$}}
\end{picture}
\end{centering}
\caption{Frequency response for different values of the forcing amplitude as indicated in the figure, $c_1=0.1$, (a) $c_3=-0.72$ and (b) $c_3=-0.8$ (corresponding to the red dots in Fig.~\ref{region_c1c3}). Dashed lines indicate unstable motions.}\label{freq_resp}
\end{figure}
For $c_3<-0.667$ the system is in zone 3. Because of the existence of several singularities, different scenarios are possible.
The system presents two IRCs already for infinitesimal values of $f$, since two isola singularities occur for $f=0$ (region d of Fig.~\ref{region_c1c3}b), depicting a frequency response curve described as ``island chain" in \cite{hirai1978general}.
The lower IRC is fully unstable, while the upper one is partially stable (see Fig.~\ref{freq_resp}).
Increasing the forcing amplitude, if $-0.746<c_3<-0.667$, the system encounters a simple bifurcation singularity which corresponds to the merging of the two IRCs, generating what we call a ``superisola".
On the contrary, for $c_3<-0.746$, the two IRCs merge one at a time with the main branch.
These scenarios are clearly shown in Fig.~\ref{freq_resp}.
The intersection of the branches of hysteresis and simple bifurcation singularities contributes to further differentiate the global picture.

The frequency response functions of Fig.~\ref{freq_resp}, obtained combining a path-following algorithm with shooting and pseudo-archlength continuation, confirms the accuracy of the analytical calculation, both qualitatively and quantitatively.
Furthermore, we stress how such non-trivial dynamical scenario was predicted with an extremely simple analytical procedure.

In all the cases studied, the frequency responses present unstable branches due to Neimark-Sacker bifurcations, therefore, the system presents quasiperiodic branches, the analysis of which is out of the scope of this paper.

\section{Geometrical interpretation}

The overall scenario illustrated in the previous sections suggests two main mechanisms leading to IRCs.
In Fig.~\ref{numerical_results}a the IRC appears for a specific forcing amplitude, while in Fig.~\ref{freq_resp} the two IRCs are generated already for infinitesimal values of $f$.

In the latter case, the leading mechanism is related to the negative damping force for a specific velocity range. Due to this, damping introduces energy in the system instead of dissipating it, causing large amplitude oscillations.
This phenomenon resembles self-excited oscillations, although it is still related to a resonance between an external excitation and the system.
In the former case, the non-monotonically increasing shape of the damping force function seems to be a critical factor for IRCs generation.

\begin{figure}
\begin{centering}
\setlength{\unitlength}{\textwidth}
\begin{picture}(1,0.35)
\put(0.0,-0.01){\includegraphics[trim = 8mm 10mm 19mm 10mm,clip,width=0.32\textwidth]{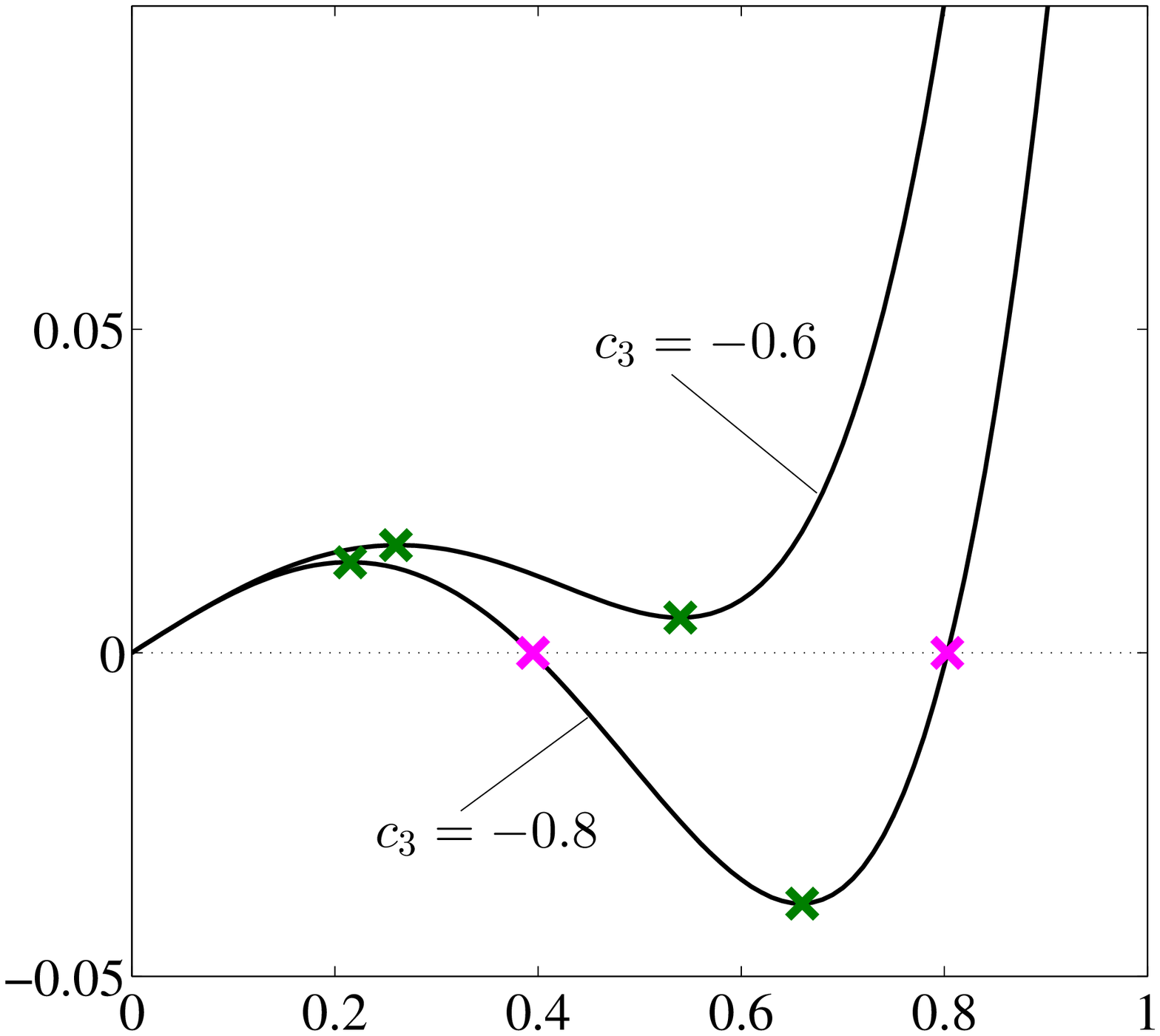}}
\put(0.35,-0.01){\includegraphics[trim = 15mm 10mm 12mm 10mm,clip,width=0.32\textwidth]{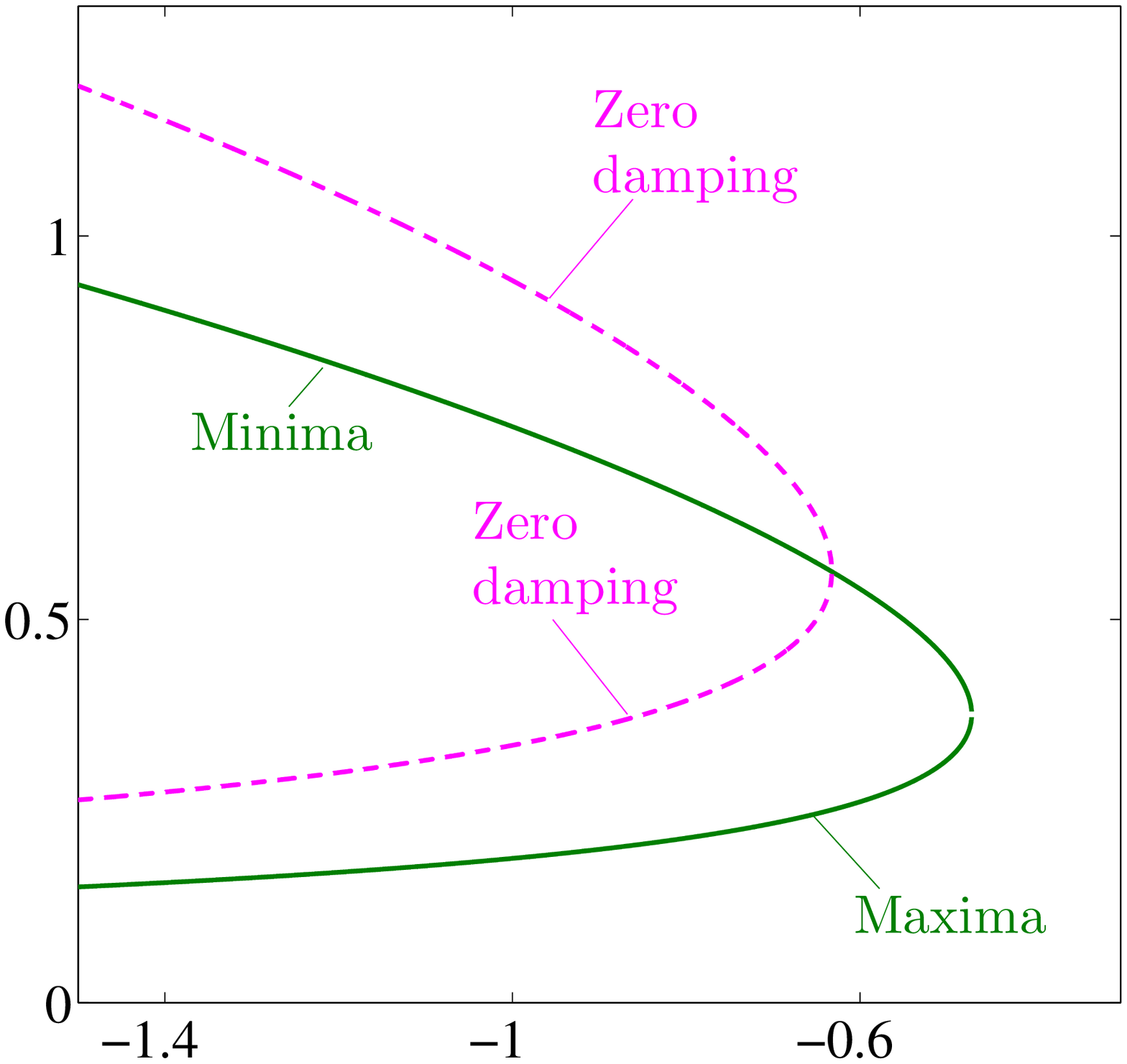}}
\put(0.69,-0.01){\includegraphics[trim = 15mm 10mm 12mm 10mm,clip,width=0.32\textwidth]{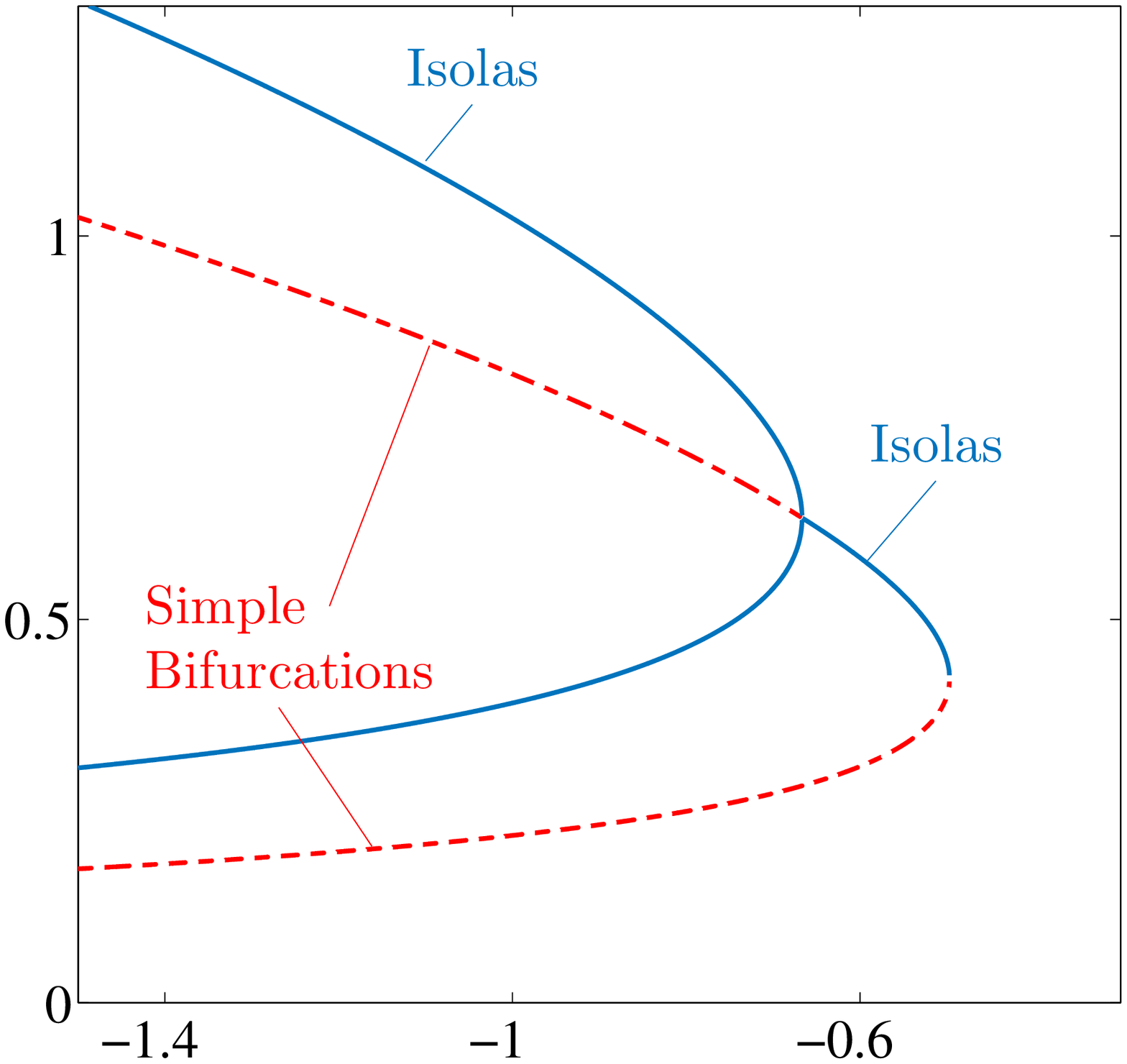}}
\put(0.04,0.26){(a)}
\put(0.38,0.26){(b)}
\put(0.72,0.26){(c)}
\put(0.175,-0.03){$\dot x$}
\put(0.51,-0.03){$c_3$}
\put(0.85,-0.03){$c_3$}
\put(-0.01,0.14){\rotatebox{90}{$F_d$}}
\put(0.33,0.14){\rotatebox{90}{$\dot x$}}
\put(0.67,0.14){\rotatebox{90}{$\dot x$}}
\end{picture}
\end{centering}
\caption{(a) Damping force for the system in Eq.~(\ref{NL_DF}); (b) maxima and minima (green solid lines) and points of zero damping (magenta dashed lines) of the damping force; (c) isolas (solid blue lines) and simple bifurcations singularities (dashed red lines) in the $c_3,\dot x$ space. For all plots $c_1=0.1$.}\label{fig_damp_force}
\end{figure}

This speculative argumentation hints that a simple graphical analysis of the damping force might already give some insight into the dynamics of the system.
Fig.~\ref{fig_damp_force}a depicts the damping force for parameter values considered in Figs.~\ref{numerical_results}a and \ref{freq_resp}b.
Maxima, minima and zeros of the curves (excluding the trivial one) are marked by green and magenta crosses.
In Fig.~\ref{fig_damp_force}b maxima and minima (solid green lines) and zeros (dashed magenta lines) of the damping force are plotted with respect to $c_3$.
In Fig.~\ref{fig_damp_force}c, isola and simple bifurcation singularities (solid blue lines and red dashed lines, respectively), obtained through single-harmonic approximation, are represented.
The qualitative similarity of the curves in Fig.~\ref{fig_damp_force}b,c suggests a strong connection between the two.
In particular, if the damping force does not cross zero for non-zero velocities (i.e., for $c_3>-0.63$ in Fig.~\ref{fig_damp_force}b), minima and maxima of the damping curve correspond to isola singularities and simple bifurcations, respectively. If there are one or more zero crossings, they coincide with isola singularities and the minimum between them is a simple bifurcation.

In spite of the good qualitative agreement, there is a non-negligible mismatch between Figs.~\ref{fig_damp_force}b and \ref{fig_damp_force}c.
The discrepancy between the two is quantified in Tab.~\ref{Compare} for the case of $c_1=0.1$ and $c_3=-0.8$. The first two columns of the table refer to numerical and analytical (through singularity theory) estimation of onset and merging points of IRCs. The matching between the two is excellent, with differences of the order of 0.05~\%. The third columns indicates zeros, maxima and minima of the damping force.
In this case the mismatch is significant, being of the order of 15~\%.
This difference is due to the fact that maxima and minima of the damping force function are referred to a specific velocity, while the periodic solution existing at a singularity has a periodically varying velocity and the singularity point indicates only its maximum.

\begin{table}
\centering
\begin{tabular}{ccccc}
&numerical&singularities&damping&av. damping\\
&$\left(x,2f\right)$&$\left(x,2f\right)$&$\left(\dot x,F_d\right)$&$\left(x_0,f_d\right)$\\
\hline
Onset 1 / zero 1&(0.4635,0) &(0.4633,0)&(0.3938,0)& (0.4633,0)\\
Onset 2 / zero 2&(0.8652,0) &(0.8633,0)&(0.8031,0) &(0.8633,0)\\
Merge 1 / maximum &(0.2497,0.01623)&(0.2496,0.01624)& (0.2147,0.01401)&(0.2496,0.01624)\\
Merge 2 / minimum& (0.7169,0.03106)&(0.7167,0.03103) &(0.6587,0.03877) &(0.7167,0.03103)
\end{tabular}
\caption{Points of onset and merging of IRCs (zeros, maximum and minimum of damping and averaged damping curves) for $c_1=0.1$ and $c_3=-0.8$ according to numerical computation and analytical prediction through singularity identification.}
\label{Compare}
\end{table}


In order to take into account this effect, instead of computing directly the damping force, we calculate the forcing amplitude necessary to compensate damping at a given oscillation amplitude. We take advantage of the energy balance criterion \cite{hill2014analytical}, which states that at resonance, for a periodic response, the energy dissipated by the system damping is equal over one period to the energy introduced into the system by the excitation force, i.e. 
\begin{equation}
\int_0^T\dot x(t)F_d\left(\dot x(t)\right)\text dt=\int_0^T\dot x(t)F\left(t\right)\text dt \label{EnBal},
\end{equation}
where $F_d$ is the damping force and $F\left(t\right)$ is the excitation force.
Continuing to approximate the motion as mono-harmonic, we impose $x=x_0\sin(\omega t)$ (with $\omega=1$ since elastic force is linear and resonance is approximately at 1).
We obtain
\begin{equation}
\int_0^T x_0\cos\left(t\right)\left(c_1x_0\cos\left(t\right)+c_3\left(x_0\cos\left(t\right)\right)^3+\left(x_0\cos\left(t\right)\right)^5\right)\text dt=\int_0^T x_0\cos\left(t\right)2f\cos\left(t\right)\text dt \label{EnBal2},
\end{equation}
computing the integrals, we have \begin{equation}
f_d=2f=c_1x_0+\frac{3}{4}c_3x_0^3+\frac{5}{8}x_0^5.\label{fx0}
\end{equation}
Eq.~(\ref{fx0}) indicates (approximately), for a given oscillation amplitude $x_0$, the forcing amplitude required to compensate damping force. We call this quantity averaged damping $f_d$.

The forth column of Tab.~\ref{Compare} indicates the zeros, maximum and minimum points of the averaged damping curve for $c_1=0.1$ and $c_3=-0.8$.
These correspond exactly to the estimated values of the singularities and approximate very well the appearance and merging of the IRCs.
This excellent matching proves the relation between zeros and extrema of the damping curve and onset and merging of IRCs.
We notice that the method adopted in \cite{kuether2015nonlinear, hill2016analytical}, which also exploits the energy balance procedure, would provide the same result.

%
\begin{figure}
\begin{centering}
\setlength{\unitlength}{\textwidth}
\begin{picture}(1,0.5)
\put(0.28,-0.01){\includegraphics[trim = 8mm 10mm 14mm 10mm,clip,width=0.45\textwidth]{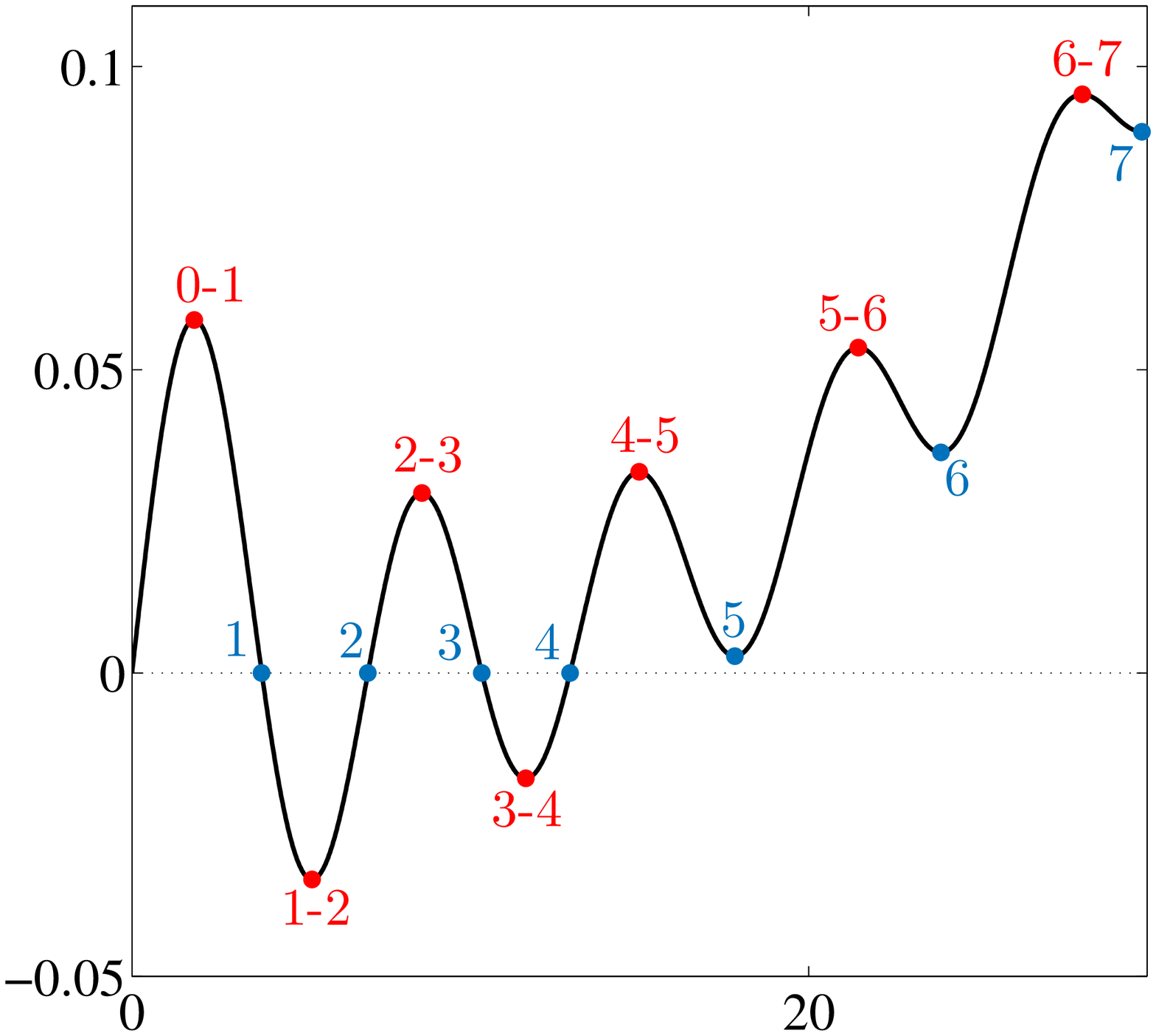}}
\put(0.5,-0.03){$x_0$}
\put(0.25,0.2){\rotatebox{90}{$f_d$}}
\end{picture}
\end{centering}
\caption{Filtered damping force for the system in Eq.~(\ref{eq_example}) for $c_1=0.1$ and $c_3=10^{-5}$.}\label{example}
\end{figure}
\begin{figure}
\begin{centering}
\setlength{\unitlength}{\textwidth}
\begin{picture}(1,0.5)
\put(0.22,-0.01){\includegraphics[trim = 10mm 10mm 12mm 10mm,clip,width=0.55\textwidth]{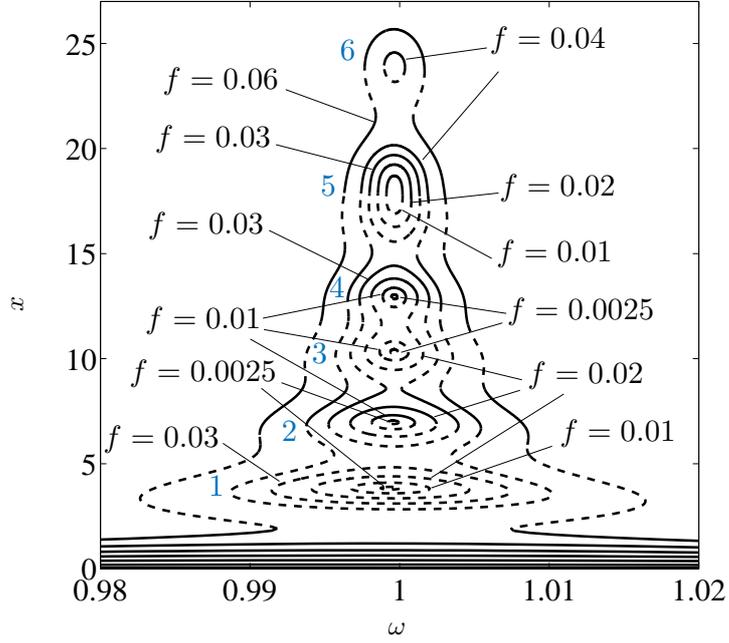}}
\put(0.5,-0.03){$\omega$}
\put(0.2,0.23){\rotatebox{90}{$x$}}
\end{picture}
\end{centering}
\caption{Frequency response of the system in Eq.~(\ref{eq_example}) for $c_1=0.1$ and $c_3=10^{-5}$.}\label{Tree}
\end{figure}

In order to validate this geometrical interpretation, we consider a system possessing a damping force having several zeros, maxima and minima, namely
\begin{equation}
\ddot x+x+c_1\sin\left(\dot x\right)+c_3\dot x^3=2f\cos\left(\omega t\right).\label{eq_example}
\end{equation}
Adopting the energy balance criterion and a single harmonic approximation, results in the averaged damping force \begin{equation}
f_d=c_1J_1\left(x_0\right)+\frac{3}{8}c_3x_0^3,\label{bessel}
\end{equation}
where $J_1$ is the Bessel function of the first kind.
Eq.~(\ref{bessel}) is depicted in Fig.~\ref{example}.
According to our geometrical interpretation, we can deduce from Fig.~\ref{example}a the appearance and merging of the IRCs. These are summarized in Tab. \ref{Merge} for $c_1=0.1$ and $c_3=10^{-5}$.
\begin{table}
\centering
\begin{tabular}{ccc}
\multicolumn{3}{c}{Onset}\\
\hline
$\#$&$x$&$f$\\
\hline
1&3.83&0\\
2&6.97&0\\
3&10.34&0\\
4&12.95&0\\
5&17.82&0.00278\\
6&23.91&0.0364\\
7&29.85&0.08923\\
...&...&...
\end{tabular}
\begin{tabular}{ccc}
\multicolumn{3}{c}{Merging}\\
\hline
$\#-\#$&$x$&$f$\\
\hline
3-4&11.64&0.0174\\
2-3&8.57&0.0297\\
4-5&14.99&0.0332\\
1-2&5.32&0.0341\\
5-6&21.47&0.0537\\
0-1&1.84&0.0582\\
6-7&28.09&0.0954\\
...&...&...
\end{tabular}
\caption{Points of onset and merging of IRCs of system in Eq.~(\ref{eq_example}) for $c_1=0.1$ and $c_3=10^{-5}$. $\#-\#$ indicates between which IRCs the merging occur, 0 refers to the main branch.}
\label{Merge}
\end{table}
The corresponding frequency response is depicted in Fig.~\ref{Tree} for increasing values of forcing amplitude $f$.

For arbitrarily small forcing amplitude, we expect to have four IRCs, because of the four zeros of Eq.~(\ref{bessel}) for $x_0>0$ (marked in Figs.~\ref{example} by the blue numbers 1, 2, 3 and 4).
This is confirmed numerically in Fig.~\ref{Tree}, where, for $f=0.0025$, four IRCs are identified.
Slightly increasing the forcing amplitude such that $f>0.00278$, the IRC number 5 is generated by the corresponding minimum of the averaged damping curve; this is verified numerically for $f=0.01$.
For $f=0.0174$, the minimum between the IRCs number 3 and 4 is reached, making the two merge, as illustrated in Fig.~\ref{Tree} for $f=0.02$.
For $f=0.0297$, also IRC number 2 merges with the 3 and 4, as shown numerically for $f=0.03$.
A further increase of $f$ up to 0.04, causes the merging between IRCs number 4 and 5 and of IRCs number 1 and 2, while for $f=0.0364$ IRC number 6 is generated.
Finally, reaching $f=0.06$, all the IRCs considered so far merge with each other, resulting in a connected frequency response curve.
Increasing even more the forcing amplitude, other minima of the averaged damping are encountered, causing the onset of new IRCs (as IRC number 7 for $f=0.08923$), which are not illustrated in Fig.~\ref{Tree}.
We point out that also the amplitude of appearance and merging of the IRCs is accurately predicted by zero, maximum and minimum points of the averaged damping curve.

Despite the very unusual shape of the frequency response, this simple procedure very precisely predicts the onset and merging of the IRCs, providing important insight into the dynamics of the system.
We also notice that the shape of the filtered damping force seems to provide information about stability as well. In fact, IRCs generated by a zero for decreasing damping are fully unstable, all the others have their upper branch stable.

\section{Conclusions}

Our purpose in this study was to understand the relationship between nonlinear damping and the existence of IRCs.
Beside several known mechanisms leading to the generation of IRCs, such as internal resonances, symmetry breaking, discontinuities and friction, it turns out that nonlinear damping is able to generate IRCs when the associated restoring force presents extrema.
Exploiting singularity theory, we were able to clearly identify parameter values for which the generation of IRCs is possible.
The uncommon scenario of two IRCs merging with each other, forming a so-called ``superisola", was also unveiled.

The obtained results illustrated that a strong connection between the shape of the damping force and IRCs onset exists.
Exploiting this relationship, we were able to predict, with a simple procedure, appearance and merging of several IRCs in system encompassing a complicated damping force function, showing that the identified connection is a general rule for this kind of systems.

\section*{Acknowledgements}

This study was financially supported by the European Union, H2020 Marie Sk\l{}odowska-Curie IF 704133 (G. Habib) and ERC Starting Grant NoVib 307265 (G. Kerschen).

\bibliographystyle{ieeetr}
\bibliography{references}

\end{document}